\documentclass[11pt]{article}
\usepackage{graphicx}
\usepackage{amsmath,amssymb,amsthm,amsfonts}
\usepackage{color}
\usepackage{cite}
\newtheorem{thm}{Theorem}[section]
\newtheorem{cor}[thm]{Corollary}
\newtheorem{lem}[thm]{Lemma}

\theoremstyle{definition}

\theoremstyle{remark}
\usepackage{appendix}

\numberwithin{equation}{section}

\DeclareMathSymbol{\C}{\mathalpha}{AMSb}{"43}

\newcommand{\eps}{\varepsilon}

\newcommand{\alp}{\alpha}
\newcommand{\R}{{\mathbb{R}}}
\newcommand{\h}{{\mathcal{H}}}
\newcommand{\inte}{\int_{\mathbb{R}^2}}

\newcommand{\e}{\varepsilon_n}

\newcommand{\bbint}{{-\mkern -16mu\int}}

\def\R{{\mathbb R}}
\def\C{{\mathbb C}}

\newcommand{\bsub}{\begin{subequations}}
\newcommand{\esub}{\end{subequations}$\!$}

\begin{document}
\title{Existence and Asymptotic Behavior of Ground States for Rotating Bose-Einstein Condensates}
\author{Yujin Guo\thanks{School of Mathematics and Statistics, and Hubei Key Laboratory of
Mathematical Sciences, Central China Normal University, P.O. Box 71010, Wuhan 430079, P. R. China. Email: \texttt{yguo@mail.ccnu.edu.cn}. Y. Guo is partially supported by NSFC under Grant No. 11931012 and the Fundamental Research Funds for the Central Universities (No. KJ02072020-0319).
},
\, Yong Luo\thanks{School of Mathematics and Statistics, and Hubei Key Laboratory of
Mathematical Sciences, Central China Normal University, P.O. Box 71010, Wuhan 430079,
P. R. China.  Email: \texttt{yluo@mail.ccnu.edu.cn}. Y. Luo is partially supported by the Project funded by China Postdoctoral Science Foundation No. 2019M662680.}
\, and\, Shuangjie Peng\thanks{School of Mathematics and Statistics, and Hubei Key Laboratory of
Mathematical Sciences, Central China Normal University, P.O. Box 71010, Wuhan 430079,
P. R. China.  Email: \texttt{sjpeng@mail.ccnu.edu.cn}. S. Peng is partially supported by the Key Project of NSFC under Grant No. 11831009.
}
}

\date{\today}

\smallbreak \maketitle
\begin{abstract}
We study ground states of two-dimensional Bose-Einstein condensates with repulsive ($a>0$) or attractive ($a<0$) interactions in a trap $V (x)$ rotating at the velocity $\Omega $. It is known that there exist critical parameters $a^*>0$ and $\Omega ^*:=\Omega^*(V(x))>0$ such that if $\Omega>\Omega^*$, then there is no ground state for any $a\in\R$; if $0\le \Omega <\Omega ^*$, then ground states exist if and only if $a\in(-a^*,+\infty)$. As a completion of the existing results, in this paper, we focus on the critical case where $0<\Omega=\Omega^*<+\infty$ and classify the existence and nonexistence of ground states for $a\in\R$. Moreover, for a suitable class of radially symmetric traps $V(x)$, employing  the inductive symmetry method, we prove that up to a constant phase, the ground states must be real-valued, unique and free of vortices as $\Omega \searrow 0$, no matter whether the interactions of the condensates are repulsive or not.
\end{abstract}	

\vskip 0.05truein

%\noindent {\it  MSC}: 35J60; 35Q40; 46N50


\vskip 0.2truein

\section{Introduction}
%Bose-Einstein condensates (BECs) is a rare state (or phase) of matter in which a large percentage of bosons collapse into their lowest quantum state.
%Bose-Einstein Condensate (BECs) is a state of matter that occurs when a large fraction of atoms suddenly occupy the lowest energy quantum state below a critical temperatures near to absolute zero.
Because Bose-Einstein condensates (BECs) can show the quantum effects at the macroscopic scale, they have become an important subject in experimental investigations since their first realization in dilute gases of alkali atoms in 1995 \cite{Abo,A,Lewin,D}.
Various interesting quantum phenomena have been so far observed in BECs experiments, such as the critical-mass collapse \cite{D}, the appearance of quantized vortices \cite{Abo}, the center-of-mass rotation in rotating traps \cite{Abo,LC}, and so on. In particular, the complex structures of BECs in rotating traps have been observed and analyzed extensively  (cf. \cite{Abo,A,AN,CRY,D,F}) starting from the late 1990s.

%An important issue in BECs theory is the relationship between BECs and superfluidity, in particular through the existence of vortices in a rotating trap \cite{A,F}.
As addressed in \cite{F,Lewin,LSY}, the two-dimensional BECs in  a rotating trap $V(x)$ can be well described
by the complex-valued wave functions minimizing the following Gross-Pitaevskii (GP) energy functional:
\begin{equation}
F_{\Omega,a}(u):=\int _{\R ^2} \big(|\nabla u|^2+V(x)|u|^2+\frac{a}{2}|u|^4\big)dx-
\Omega \int_{\R ^2}x^{\perp}\cdot (iu,\, \nabla u)dx,  \ \ u\in \h , \label{f}
\end{equation}
where the space $\h$ is defined as
\begin{equation}\label{1.H}
\h :=  \Big \{u\in  H^1(\R ^2, \mathbb{C}):\ \int _{\R ^2}
V(x)|u|^2 dx<\infty\Big \},
\end{equation}
$x^{\perp} =(-x_2,x_1)$ with $x=(x_1,x_2)\in \R^2$, and
$(iu,\, \nabla u)=i(u\nabla \bar u-\bar u\nabla u)/2$.
Here the two-body interaction of the rotating BECs can be either repulsive ($a>0$) or  attractive ($a<0$), and the parameter $\Omega \ge 0$ describes the rotational velocity of the rotating trap. Consequently, the {\bf ground states} of two-dimensional  BECs in a rotating trap $0\leq V (x)\in L^{\infty}_{loc}(\R^2)$ can be equivalently described (see \cite{F,Lewin,LSY,P}) by the {\bf minimizers} of the following variational problem:
\begin{equation}\label{def:ea}
e(\Omega,a):=\inf _{\{u\in \h, \, \|u\|^2_2=1 \} } F_{\Omega,a}(u),\,\ \hbox{$\Omega\geq 0$, $a\in\R$},
\end{equation}
where the GP energy functional $F_{\Omega,a}(\cdot)$ is as in (\ref{f}), and the space $\h$ is defined by (\ref{1.H}). Here the value $|a|$ characterizes the absolute product for the scattering length  of the two-body interaction times the number $N$ of particles in the condensates.

We first introduce the non-rotational case $\Omega=0$ of $e(\Omega,a)$. As illustrated soon after (\ref{2.27}) (see also \cite[Section 2]{GLY}), in this case it suffices to consider the real-valued minimizers of $e(a):=e(0,a)$, where we define
\begin{equation}\label{2.2}
e(a):=\inf _{\{u\in \h, \, \|u\|^2_2=1 \} } F_{a}(u),
\end{equation}
and the GP energy functional $F_{a}(u)$ is given by
\begin{equation}\label{2.1}
F_{a}(u):=\inte \big(|\nabla u|^2+V(x)|u|^2+\frac {a}{2} |u|^4\big)dx.
\end{equation}
When $a>0$, since the functional $F_a(u)$ is convex, one can obtain the existence and uniqueness of  real-valued minimizers for $e(a)$ in a direct way, see for instance \cite[Theorem 2.1]{LSY}. However, when $a<0$, it was shown in \cite{GLW,GS,GWZZ} that there exists a critical constant $a^*>0$ such that $e(a)$ admits real-valued minimizers, if and only if $a>-a^*$. More precisely, the critical constant $a^*$ is given by
\begin{equation}\label{1.3}
a^*:=\|w\|^2_{L^2(\R^2)}>0,
\end{equation}
where $w=w(|x|)>0$ is a positive solution of the following scalar field equation
\begin{equation}
-\Delta w+ w-w^3=0\  \mbox{  in } \  \R^2,\  \mbox{ where }\ w\in H^1(\R ^2,\R), \label{Kwong}
\end{equation}
which, up to translations, must be unique and radially symmetric (cf. \cite{K,W}).
Moreover, the refined mass concentration and local uniqueness of real-valued minimizers of $e(a)$ as $a\searrow -a^*$ were also studied in \cite{GLW,GWZZ} and the references therein.

%More recently, the variational problem $e_F(a)$ under the rotation $\Omega >0$ was discussed in \cite{BC,ANS,GLY,Lewin}, where the authors studied the existence, non-existence and the limit behavior of complex-valued minimizers.

As for the rotational case $\Omega>0$ of $e(\Omega,a)$, we consider the general trapping potential $0\leq V (x)\in L^{\infty}_{loc}(\R^2)$ satisfying
\begin{equation}\label{A:V}
\underline{\lim} _{|x|\to\infty }\frac{V(x)}{|x|^2}>0.
\end{equation}
Under the assumption (\ref{A:V}), we define the critical rotational velocity $\Omega ^*=\Omega ^*(V(x))$ by
\begin{equation}
\Omega ^*:=\sup \Big\{\Omega >0:\ \  V(x)-\frac{\Omega ^2}{4}|x|^2  \to\infty \,\ \mbox{as}\,\
|x|\to\infty \Big\},  \label{Omega}
\end{equation}
so that $\Omega ^*>0$.
Considering $V(x)=|x|^s$ $(s\ge 2)$  as an example, we have
\begin{equation}
\Omega ^* :=\arraycolsep=1.5pt \left\{ \begin{array}{ll}
2,\quad &  {\rm if}\quad s=2;\\ [2mm]
\infty, \ &  {\rm if}\quad s>2,\end{array}\right.
\label{omega-1}
\end{equation}
which illustrates that  both $0<\Omega ^*<\infty$ and $\Omega ^*=\infty$ can happen. Recalling from \cite{GLY,Lewin,S02}, we have the following existence and non-existence of minimizers:

\vskip 0.05truein
\noindent{\bf Theorem A. } {\em  (\cite{GLY,Lewin,S02}) Assume $V(x)\in L^\infty_{\rm loc}(\R^2)$ satisfies
	(\ref{A:V}), and let $\Omega^*>0$ and $a^*>0$ be defined in \eqref{Omega} and (\ref{1.3}), respectively. Then we have
\begin{enumerate}
\item If $\Omega\in(0,\Omega ^*)$ and $a\in(-a^*,+\infty)$, then there exists at least one minimizer for $e(\Omega,a)$.
\item If $\Omega\in(0,\Omega ^*)$ and $a\in(-\infty,-a^*]$, then there is no minimizer for $e(\Omega,a)$.
\item If $\Omega>\Omega^*$, then there is no minimizer of $e(\Omega,a)$ for any $a\in\R$.
\end{enumerate}}

%We remark that if $\Omega>\Omega^*$, then one can prove as in \cite{GLY,Lewin} that $e(\Omega,a)=-\infty$ holds for any $a\in\R$,
%which thus gives Theorem A (3). However, the case $a\ge 0$ of Theorem A (1) can be established in the similar way of \cite{IM-1},
%while the rest cases of Theorem A can be proved by borrowing the ideas of \cite{GLY,Lewin}. For simplicity we omit the detailed
%proof of Theorem A.

One can note from Theorem A that it remains open to discuss the existence and non-existence of minimizers for $e(\Omega ^*,a)$, provided that the critical rotational velocity $\Omega^*=\Omega ^*(V(x))$ defined in (\ref{Omega}) satisfies $\Omega^*<+\infty$. The first main purpose of this paper is to address this issue.

When $0<\Omega<\Omega^*$, the detailed analytical properties of minimizers for $e(\Omega,a)$ as either $a\nearrow+\infty$ or $a\searrow-a^*$  were studied recently in \cite{A,IM-1,IM-2,CRY,GLP,GLY,Lewin} and the references therein. Specially, if the interactions in the condensates are repulsive ($a>0$), all kinds of  quantized vortices for $e(\Omega,a)$ have been analyzed extensively in the Thomas-Fermi regime where $a\nearrow+\infty$, see \cite{Abo,A,AJ,AN,CRY,D,IM-1,IM-2} for instance. In particular, the nonexistence of  vortices for repulsive BECs under rotation were analyzed in \cite{A,AJ,AN}. However, as far as we know, the above mentioned works of studying vortex structures focus more on the analysis of the energy, which seems not enough for the attractive case $a<0$ of $e(\Omega,a)$. Fortunately, the nonexistence of vortices for $e(\Omega,a)$ in the attractive limiting case as $a\searrow-a^*$ was proved in our recent work \cite{GLY} by developing the inductive symmetry method. The second main purpose of this paper is to address the nonexistence of  vortices for $e(\Omega,a)$ as $\Omega \searrow 0$ by a uniform approach, no matter whether the fixed parameter $a>-a^*$ is positive or not.

\subsection{Main results}
In this subsection we shall introduce the main results of the present paper. We are first concerned with the existence and non-existence of minimizers for $e(\Omega ^*,a)$, provided that the critical rotational velocity $\Omega^*=\Omega ^*(V(x))$ defined in (\ref{Omega}) satisfies $\Omega^*<+\infty$.
Towards this aim, we consider
\begin{equation}\label{1.2}
V(x)=A|x|^2+W(x),\,\ \hbox{where $A>0$ and $W(x)\in C(\R^2)$ satisfies $\lim_{|x|\to \infty}\frac{W(x)}{|x|^2}=0$,}
\end{equation}
so that the assumption (\ref{A:V}) holds in a natural way.
Under the assumption (\ref{1.2}), our first main result of this paper is devoted to the following existence and non-existence of minimizers for $e(\Omega^*,a)$:

\begin{thm}\label{thm1.1}
Suppose that $a^*>0$ is defined by \eqref{1.3}, and $V(x)$ satisfies \eqref{1.2} such that $\Omega^*=2\sqrt{A}>0$ is defined by \eqref{Omega}. Then we have
\begin{enumerate}
\item[\rm(i).] If $W(x)$ satisfies $\lim_{|x|\to\infty }W(x)=+\infty$, then $e(\Omega^*,a)$ admits minimizers if and only if $a\in(-a^*,\infty)$.
\item[\rm(ii).] If  $W(x)\equiv C\in\R$, then $e(\Omega^*,a)$ admits  minimizers if and only if $a\in(-a^*,0]$.
\end{enumerate}	
\end{thm}

When the non-constant function $W(x)$ is bounded uniformly in $\R^2$, the existence and non-existence of minimizers for $e(\Omega^*,a)$ are more complicated.  The following theorem addresses this case:

\begin{thm}\label{thm1.1*}
Suppose that $a^*>0$ is defined by \eqref{1.3}, and $V(x)$ satisfies \eqref{1.2} such that $\Omega^*=2\sqrt{A}>0$ is defined by \eqref{Omega}. If $W(x)\not\equiv const.$ further satisfies $\lim_{|x|\to\infty }W(x)=B\in\R$, then we have
\begin{enumerate}
\item[\rm(i).] If $W(x)>B$ in $\R^2$, then $e(\Omega^*,a)$ admits no minimizer for any $a\in \R$.
\item[\rm(ii).] If $W(x)<B$ in $\R^2$, then there exists a constant $a_\ast\in(0,+\infty]$ such that
\begin{enumerate}
\item[(a).] If $a\in(-a^\ast,a_\ast)$, there is at least one minimizer for $e(\Omega^*,a)$.
\item[(b).] If $a\in(-\infty,-a^*]\cup (a_\ast,+\infty)$, there is no minimizer for $e(\Omega^*,a)$.
\end{enumerate}
\end{enumerate}
\end{thm}

%and $\liminf_{|x|\to\infty}\big[B-W(x)\big]|x|^s>0$ ($0<s<2$), then
%$e(\Omega^*,a)$ admits at least one minimizer if and only if $a\in (-a^*,+\infty)$.
%\item[\rm(iii).]

The proof of Theorem \ref{thm1.1*} (ii) shows that the constant $a_\ast\in(0,+\infty]$ of Theorem \ref{thm1.1*} (ii) can be characterized as
\begin{equation}\label{ast}
a_\ast:=\sup\big\{a\in\R:\, e(\Omega^*,a)<2\sqrt{A}+B\,\big\}.
\end{equation}
We remark that the energy $2\sqrt{A}+B$ is maximal to guarantee that vanishing cannot occur for the minimizing sequences of $e(\Omega^*,a)$. 
Theorem \ref{thm1.1*}(ii) does not address the existence of minimizers for $e(\Omega^*,a_\ast)$, which may depend on  the shape of  $W(x)<B$ in $\R^2$.
Additionally, if $W(x)$ satisfies
\begin{equation}\label{new2}
\liminf_{|x|\to\infty}\big[B-W(x)\big]|x|^s>0,\quad\hbox{where $0<s<2$,}
\end{equation}
then the following corollary shows that $a_\ast=+\infty$  holds true.

\begin{cor}\label{cor-1.3}
Under the assumptions of  Theorem \ref{thm1.1*} (ii), if $W(x)$ satisfies \eqref{new2}, then $e(\Omega^*,a)$ admits minimizers if and only if $a\in (-a^*,+\infty)$.
\end{cor}

\noindent
Related to Corollary \ref{cor-1.3}, we expect that if (\ref{new2}) holds for $s>2$, then it might hold that $a_\ast<\infty$, which however we cannot address  rigorously at present.

As far as we know, most existing results focus mainly on the noncritical case $e(\Omega ,a)$, where $\Omega\not =\Omega^*$ and the energy functional is coercive. But  the existence and non-existence results of minimizers for the critical case $e(\Omega ^*,a)$ seem very few, no matter whether the interactions among the BEC system are repulsive or not.
The proof of above results shows actually that the analysis of the critical case $e(\Omega^*,a)$ is more challenging, compared with those involved in \cite{LSY,A,GLY} for the noncritical case $e(\Omega,a)$ where $ \Omega\not= \Omega ^*$.
The above results also illustrate  that whether the minimizers of $e(\Omega ^*,a)$ exist or not depends  subtly on the shape of the trap $V(x)$ and the parameter $a\in\R$ as well.

It should also be mentioned that the mass-subcritical version of $e(\Omega,a)$, where the nonlinear term $|u|^4$ is replaced by $|u|^p$ for $2<p<4$, was studied in the pioneering work of Esteban-Lions \cite{EL} by the concentration-compactness lemma. However, those methods used in \cite{EL} are not sufficient to prove Theorem \ref{thm1.1}, especially in the repulsive case where $a>0$. To overcome this difficulty, we shall borrow some ideas from \cite{A,AB}, where  the properties of the first eigenpairs of
\begin{equation}\label{s-2.17}
-\Delta \psi+2i\,(x^{\perp}\cdot \nabla \psi)+|x|^2\psi=\lambda\psi\quad\hbox{in}\,\ \R^2
\end{equation}
were fully employed.
Note from \cite[Theorem 2]{LP} that the first eigenvalue of \eqref{s-2.17}  is $\lambda_1=2$
and all the first eigenfunctions are given by
\begin{equation}\label{s-2.18}
S:=\Big\{e^{-\frac{|x|^2}{2}}f(x): f(x)\,\ \hbox{is any entire function
	such that}\,\ e^{-\frac{|x|^2}{2}}f(x)\in L^2(\R^2,\mathbb{C})\Big\},
\end{equation}
where the entire function means the complex analytic function on
the entire plane. In Section 2 below, we shall seek for appropriate test functions in $S$ to derive the desired energy estimates,
which then help us obtain the existence and nonexistence of minimizers for the repulsive case where $a>0$.
%As remarked in \cite{GLY}, our Theorem \ref{thm1.1} also show that whether minimimizers of $e(\Omega,a)$ exist or not in $\Omega=\Omega^*$ depending heavily on the shape of the trapping potential $V(x)$ and the parameter $a\in\R$.

We now recall the non-rotational case  of $e(\Omega,a)$ at $\Omega=0$, i.e., $e(a):=e(0,a)$ defined in (\ref{2.2}). By the variational theory, any minimizer $u_a$ of $e(a)$  satisfies the following Euler-Lagrange equation
\begin{equation}\label{2.26}
-\Delta u_a+V(x)u_a=\mu_a u_a-a|u_a|^2u_a\ \ \hbox{in}\,\ \R^2,
\end{equation}
where $\mu_a\in\R$ is the associated Lagrange multiplier satisfying
\begin{equation}\label{2.27}
\mu_a=e(a)+\frac{a}{2}\inte |u_a|^4dx.
\end{equation}
Since $F_a(u)\geq F_a(|u|)$, we obtain that $|u_a|$ is also a minimizer of $e(a)$ and $F_a(u_a)= F_a(|u_a|)$ holds. Moreover, we can deduce from (\ref{2.26}) that $|u_a|>0$ by the maximal principle. Since $\inte |\nabla u_a|^2dx=\inte \big|\nabla |u_a|\big|^2dx$ and  $|u_a|>0$, we conclude that there exists a constant $\theta\in(0,2\pi]$ such that $u_a(x)\equiv|u_a(x)|e^{i\theta}$ holds in $\R^2$.  Hence, up to a constant phase, any minimizer $u_a$ of $e(a)$ must be positive.

Under the assumption that the trap  $V(x)$ is radially symmetric and satisfies
\begin{enumerate}
	\item [\rm($V$).]
	$V(x)=V(|x|)\in C^{1,\alpha}(\R^2)$, $V'(|x|)\geq 0$, $\lim_{|x|\to\infty}V(x)=+\infty$ and there exists some constant $p\ge 2$ such that 	$|V(x)|,\,\ |\nabla V(x)|\leq C|x|^p$ as $|x|\to\infty$,
\end{enumerate}
combining \cite[Theorem 2.1]{LSY} with \cite[Corallary  1.1]{GWZZ} yields immediately the following uniqueness of positive minimizers for $e(a)$:

\vskip 0.05truein
\noindent{\bf Theorem B.} {\em (\cite{LSY,GWZZ}) Suppose that $a^*>0$ is defined by \eqref{1.3}, and $V(x)$ satisfies the assumption $(V)$. For any $a\in[0,+\infty)$ and a.e. $a\in (-a^*,0)$, then $e(a)$ admits a unique positive minimizer $u_0=u(a,0)>0$, which must be radially symmetric.
}

\vskip 0.05truein

%\noindent

Based on the uniqueness of Theorem B, the second main result of this paper is the following nonexistence of vortices for  any fixed $a\in[0,+\infty)$ and a.e. $a\in (-a^*,0)$.

\begin{thm}\label{thm1.2}
Suppose that $a^*>0$ is defined by \eqref{1.3}, $V(x)$ satisfies \eqref{A:V} and the assumption $(V)$. Then for any fixed $a\in[0,+\infty)$ or a.e. $a\in (-a^*,0)$, there exists a sufficiently small $\delta:=\delta(a)>0$ such that, up to a constant phase, all minimizers of $e(\Omega,a)$ are real-valued, unique and free of vortices for any $\Omega\in[0,\delta)$.
\end{thm}

The assumption $V'(|x|)\geq 0$ of Theorem \ref{thm1.2} can be removed in the repulsive case, since it is used only in the attractive case to guarantee the uniqueness of positive minimizers for $e(a)$.
The proof of Theorem \ref{thm1.2} gives us that for any sufficiently small $\Omega\in [0,\delta)$, any minimizer $u_{\Omega}$ of $e(\Omega,a)$ satisfies $u_{\Omega}\equiv u_0e^{i\theta}$ in $\R^2$, where $\theta:=\theta(\Omega)\in[0,2\pi)$ is a suitable constant depending on $\Omega$. Thus $|u_\Omega|=|u_0|>0$ never vanishes in $\R^2$ as $\Omega \searrow 0$, which  implies  that it is free of vortices for ground states of $e(\Omega,a)$ as $\Omega \searrow 0$.

%Note that the non-existence of vortices of minimizers for repulsive BECs ($a>0$) was obtained by Jacobian estimates \cite{AJ}, which depend more on the analysis of the energy.
We next illustrate the general strategy of proving Theorem \ref{thm1.2}. Denote $\{\Omega_n>0\}$ to be a sequence satisfying $\lim_{n\to \infty}\Omega_n=0$. Choose some suitable constant  $\theta_{n}\in[0,2\pi)$ such that the minimizer  $u_n$ of $e(\Omega_n,a)$ satisfies the following orthogonal condition
\begin{equation}\label{1.11}
\inte u_0Im (u_ne^{i\theta_n})dx\equiv 0\,\ \hbox{for all $n\in \mathbb{N}^+$,}
\end{equation}
where $u_0$ is the unique positive minimizer of $e(a)$, see Theorem B.
As the first step, we shall follow  \eqref{1.11} to prove that
\begin{equation}\label{1.12}
\lim_{n\to\infty}u_ne^{i\theta_n}=u_0\,\ \hbox{uniformly in}\,\ L^\infty(\R^2,\mathbb{C}).
\end{equation}

Rewrite now $u_ne^{i\theta_n}=q_n+ir_n$, where $q_n$ and $r_n$ denote the real and imaginary parts of $u_ne^{i\theta_n}$, respectively.
One can show that $(q_n,r_n)$ satisfies the following elliptic system:
\begin{equation}\label{1.13}
\left\{
\begin{aligned}
\mathcal{L}_n q_n&=\Omega_n (x^\bot\cdot\nabla r_n)
\ \ &\mbox{in}\,\  \R^2,\\
\mathcal{L}_n r_n&=-\Omega_n(x^\bot\cdot\nabla q_n)
\ \ &\mbox{in}\,\  \R^2,\\
\end{aligned}
\right.
\end{equation}
where the operator $\mathcal{L}_n$ is defined by
\begin{equation}\label{1.14}
\mathcal{L}_n:=-\Delta+V(x)-\mu_n+a|u_n|^2,
\end{equation}
and $ \mu_n\in \R$ is a suitable Lagrange multiplier. The refined analysis gives that $ \mu_n\to \mu_0$ as $n\to\infty$,
where $\mu_0\in\R$ is the Lagrange multiplier associated to $u_0$. By applying \eqref{1.12}, we then follow the system (\ref{1.13}) to investigate the  linearized
operators $\mathcal{L}:=-\Delta+V(x)-\mu_0+au_0^2$ and $\mathcal{N}:=-\Delta+V(x)-\mu_0+3au_0^2$ in $\R^2$. As the second step, this yields finally Corollary \ref{prop-1.1} on the first estimates of $q_n-u_0$ and $r_n$  as $n\to\infty$.

Stimulated by the inductive symmetry method (cf. \cite{GLY}),  in the third step we shall follow Corollary \ref{prop-1.1} to further prove the inductive process of Lemma \ref{lem4.5} in the following sense:  once we have an estimate
on the non-radial part of $q_n-u_0$, then a better estimate on the non-radial part of $q_n-u_0$ can
be derived. As a consequence, we finally prove that the non-radial part of $q_n-u_0$ is arbitrarily small and $r_n\equiv 0$ as  $n\to\infty$. This implies finally that $u_ne^{i\theta_n}\equiv q_n\equiv u_0$ in $\R^2$ as $\Omega_n \searrow 0$, and Theorem \ref{thm1.2} is therefore proved, see Section 4 for more details.

Similar nonexistence results of Theorem \ref{thm1.2} in the repulsive case $a>0$ were obtained in \cite{A,AJ,AN} and the references therein by applying vortex ball constructions, jacobian estimates, and some other arguments, all of which make full use of the refined energy analysis and the Ginzburg-Landau theory. As far as we know, the above mentioned arguments however seem not applicable for the attractive case of $e(\Omega , a)$, since it does not admit the variational structure of the Ginzburg-Landau type.  The above proof strategy shows that Theorem \ref{thm1.2} is established essentially by the inductive symmetry method, which was first imposed in \cite{GLY} to study the nonexistence of vortices for $e(\Omega , a)$ in a different situation where  $0<\Omega<\Omega^*$ is   fixed and $a\searrow -a^*$. This approach relies more on the refined analysis of the governing equations of minimizers. The advantage of this approach lies in the fact that Theorem \ref{thm1.2} can be proved by a uniform approach, $i.e,$  no matter whether the interactions of the condensates are repulsive ($a>0$) or not. As a byproduct, Theorem \ref{thm1.2} therefore yields  the first result that the nonexistence of vortices for repulsive BECs under rotation can be investigated by the inductive symmetry method.

This paper is organized as follows. In Section 2, we shall prove Theorem \ref{thm1.1} on the existence and non-existence of minimizers for  $e(\Omega^*,a)$.  In Section 3, we first investigate the limit behavior of $u_n$ as $\Omega_n\searrow  0$, based on which we then address the detailed analysis of the linearized problem \eqref{1.14}. In Section 4 we shall complete the proof of Theorem \ref{thm1.2} by the inductive symmetry method.

\section{Existence  of minimizers for  $e(\Omega^*,a)$}
The purpose of this section is to prove Theorem \ref{thm1.1} on the existence and nonexistence of minimizers for $e(\Omega ^*,a)$, where $\Omega^*=2\sqrt{A}>0$ is defined by \eqref{Omega} under the assumption \eqref{1.2}. Towards this purpose, we first introduce the following Gagliardo-Nirenberg inequality
\begin{equation}\label{GNineq}
\inte |u(x)|^4 dx\le \frac 2 {a^*} \inte |\nabla u(x) |^2dx \inte |u(x)|^2dx ,\
\  u \in H^1(\R ^2, \R),
\end{equation}
where the equality is attained  (cf. \cite{W}) at the unique positive radial solution $w$ of \eqref{Kwong}. Moreover, following \cite[Proposition 4.1]{GNN}, we obtain that $w=w(|x|)>0$ satisfies
\begin{equation}\label{1:id}
\inte |\nabla w |^2dx  =\inte w ^2dx=\frac{1}{2}\inte w ^4dx,
\end{equation}
and
\begin{equation}
w(x) \, , \ |\nabla w(x)| = O(|x|^{-\frac{1}{2}}e^{-|x|}) \quad
\text{as \ $|x|\to \infty$.}  \label{1:exp}
\end{equation}
Given any vector function $\mathcal{A }\in L^2_{loc}(\R^2,\R^2)$, we have the following diamagnetic inequality \cite{Lieb}:
\begin{equation}
|(\nabla -i\mathcal{A } )u|^2 \ge \big| \nabla |u|\big|^2 \,\ \hbox{a.e. on }\ \R^2,\ \, u\in H^1(\R^2,\mathbb{C}).
\label{Diam}
\end{equation}
We also recall from \cite[Lemma 2.1]{IM-1} the following compactness lemma:

\begin{lem}\label{2:lem1}  Suppose $V(x) \in L_{\rm loc}^\infty(\R^2)$ satisfies $\lim_{|x|\to \infty} V(x) = \infty$, then the embedding $\h  \hookrightarrow L^{q}(\R^2, \mathbb{C})$ is compact for any $2\le q<\infty$.
\end{lem}

Applying Lemma \ref{2:lem1}, we first prove Theorem \ref{thm1.1} as follows.

\vskip 0.05truein

\noindent{\bf Proof of Theorem \ref{thm1.1}.}  Without loss of generality, we may assume $A=1$ so that $\Omega^*=2$. If $a\in (-\infty,-a^*]$, the argument of proving \cite[Theorem 2.3]{GLY} then gives the nonexistence of minimizers for $e(2,a)=e(\Omega^*,a)$. Therefore, the rest is to consider the case where $a\in (-a^*, +\infty)$.

(i). Assume $a\in (-a^*, +\infty)$ and $\lim_{|x|\to\infty}W(x)=+\infty$, hence we get that $e(2,a)>-\infty$. By the constraint condition, we may assume $W(x)\geq 0$ for simplicity.  Let $\{u_n\}$ be a minimizing sequence of $e(2,a)$. We then deduce from \eqref{GNineq} and \eqref{Diam} that
\[
\begin{aligned}
e(2,a)&=\lim_{n\to\infty}\Big\{\inte\big(|(\nabla-ix^\perp)u_n|^2+W(x)|u_n|^2+\frac{a}{2}|u_n|^4\big)dx\\
&\geq\inte \Big(\frac{a^*+sgn(a)|a|}{2}| u_n|^4+W(x)|u_n|^2\Big)dx,
\end{aligned}
\]
which implies that $\{u_n\}$ is bounded uniformly in $L^4(\R^2,\mathbb{C})\cap H^1_A(\R^2,\mathbb{C})$, where $H^1_A(\R^2,\mathbb{C})$ is defined in (\ref{s-2.26}).

Taking a subsequence if necessary, now we may assume that $u_n\rightharpoonup u_0$ weakly in $H^1_A(\R^2,\mathbb{C})$ as $n\to\infty$. On the other hand, since $\lim_{|x|\to \infty}W(x)=+\infty$, we also obtain from Lemma \ref{2:lem1} that up to a subsequence if necessary, $\{|u_n|\}$ converges to $\{|u_0|\}$ strongly in $L^p(\R^2,\R)$ $(2\leq p<\infty)$. We thus have
$\|u_0\|_{L^2(\R^2)}^2=1$ and $\lim_{n\to\infty}\|u_n\|_{L^4(\R^2)}=\|u_0\|_{L^4(\R^2)}$.  By the weak lower semicontinuity,  we further obtain that $e(2,a)=F_{2,a}(u_0)$, and hence $e(2,a)$ admits a minimizer $u_0$ if $a\in(-a^*,+\infty)$.

\vskip 0.1truein
(\rm ii). Without loss of generality, by the constraint we may assume $W(x)\equiv C=0$ so that $\Omega^*=2$, and
\begin{equation}\label{2:newF}
F_{2,a}(u):=\int_{\mathbb{R}^2} |(\nabla-ix^{\perp})u|^2dx+\frac{a}{2}\int _{\R ^2}  |u|^4dx,  \ \ u\in \h .
\end{equation}
The proof is then divided into two cases: $a\in(-a^*,0]$ and $a\in(0,+\infty)$.

(1). We first consider the case $a\in(-a^*,0]$. The existence of minimizers of $e(2,a)=e(\Omega^*,a)$ at $a=0$ is trivial. Indeed, it follows from \cite[Theorem 2]{LP} (see also \cite[Remark 2.5]{EL}) that if $u\in\mathcal{H}$ satisfies $\|u\|^2_2=1$, then
\begin{equation}\label{s-2.19}
\int_{\mathbb{R}^2} |(\nabla-ix^{\perp})u|^2dx\geq 2\int_{\mathbb{R}^2} |u|^2dx=2,
\end{equation}
where the equality is attained at $u=\frac{1}{\sqrt{\pi}}e^{-\frac{|x|^2}{2}}$. On the other hand, by choosing $\frac{1}{\sqrt{\pi}}
e^{-\frac{|x|^2}{2}}$ as a test function, one can obtain that   $e(2,0)\le 2$, which hence yields that
$e(2,0)= 2$. Therefore, $e(2,0)$ admits at least one minimizer which is $\frac{1}{\sqrt{\pi}}e^{-\frac{|x|^2}{2}}$.
%Due to the constraint $\|u\|^2_{L^2(\R^2)}=1$, we obtain $e(2,0)\geq 2$. Thus, one can check that $e(2,0)= 2$ and $e(2,0)$ admits a minimizer by using $u=\frac{1}{\sqrt{\pi}}e^{-\frac{|x|^2}{2}}$ as a test function.

Next, we shall use the concentration-compactness lemma \cite{Lions1} to deal with the case $a\in(-a^*,0)$. Note from \eqref{GNineq} and \eqref{Diam} that for any $u\in \mathcal{H}$ satisfying $\|u\|^2_2=1$,
\begin{equation}\label{s-2.1}
F_{2,a}(u)\geq \inte \big|\nabla |u|\big|^2dx+\frac{a}{2}\inte|u|^4dx\geq \frac{a^*-|a|}{2}\inte|u|^4dx.
\end{equation}
Therefore, $e(2,a)$ is bounded from below in this case. Let $\{u_n\}$ be a minimizing sequence of $e(2,a)$, i.e., $\|u_n\|^2_2=1$ and
\[
e(2,a)=\lim_{n\to\infty}F_{2,a}(u_n)=\lim_{n\to\infty}\Big\{\inte \Big(|(\nabla -ix^{\perp})u_n|^2-\frac{|a|}{2}|u_n|^4\Big)dx\Big\}.
\]
We then derive from \eqref{s-2.1} that there exists a constant $C(a)>0$, depending only on $a$, such that $\{u_n\}$ satisfies
\begin{equation}\label{s-2.2}
\inte|u_n|^4dx\leq C(a)\,\ \hbox{and}\,\ \inte |(\nabla -ix^{\perp})u_n|^2dx\leq C(a).
\end{equation}
Denote
\begin{equation}\label{s-2.3}
f_n(x):=|(\nabla -ix^{\perp})u_n|^2+|u_n|^2\ge 0,
\end{equation}
so that the nonnegative sequence $\{f_n\}$ is bounded uniformly in $L^1(\R^2)$ in view of \eqref{s-2.2}. By the concentration-compactness lemma \cite[Lemma I.1]{Lions1}, either the vanishing or the dichotomy or the compactness holds for some subsequence of $\{f_n\}$.

We first claim that the vanishing
\begin{equation}\label{s-2.4}
\lim_{n\to\infty}\sup_{y\in \R^2}\int_{B_R(y)}f_n(x)dx= 0,\quad \forall R>0
\end{equation}
cannot occur. On the contrary, suppose that \eqref{s-2.4} holds true.
We then get from \cite[Lemma 1.21]{WM} that up to a subsequence if necessary,
\begin{equation}\label{s-2.5}
\lim_{n\to\infty}\inte|u_n|^4dx= 0.
\end{equation}
Following \cite[Theorem 2]{LP}, we thus obtain from \eqref{s-2.5} that
\[
e(2,a)=\lim_{n\to\infty}\Big\{\inte \Big(|(\nabla -ix^{\perp})u_n|^2-\frac{|a|}{2}|u_n|^4\Big)dx\Big\}\geq 2,
\]
which however contradicts to the estimate
\begin{equation}\label{s-2.27}
e(2,a)\leq F_{2,a}\Big(\frac{1}{\sqrt{\pi}}e^{-\frac{|x|^2}{2}}\Big)=2-\frac{|a|}{\pi^2}\inte e^{-2|x|^2}dx<2.
\end{equation}
Hence, the vanishing cannot occur.

To rule out the dichotomy, we first claim that
\begin{equation}\label{s-2.6}
e(2,a)< e_{\lambda}(2,a)+e_{1-\lambda}(2,a)\qquad\hbox{for any $\lambda\in (0,1)$,}
\end{equation}
where $e_{\lambda}(2,a)$ is defined by
\begin{equation}\label{s-2.7}
e_{\lambda}(2,a):=\inf _{\{u\in \h, \, \|u\|^2_2=\lambda\} } F_{2,a}(u),\quad 0<\lambda<1.
\end{equation}
Indeed, one can check from the definition of $F_{2,a}$ that
\begin{equation}\label{s-2.10}
\lambda F_{2,a}(u)=F_{2,a}(\sqrt{\lambda}u)-\frac{|a|\lambda(1-\lambda)}{2}\inte |u|^4dx\qquad\hbox{for any $\lambda\in (0,1)$.}
\end{equation}
Moreover, let $\{\sqrt{\lambda}v_{1,n}\}$ and $\{\sqrt{1-\lambda}v_{2,n}\}$ be the minimizing sequences of $e_\lambda(2,a)$ and $e_{1-\lambda}(2,a)$, respectively. The argument of proving the above vanishing gives that
\begin{equation}\label{s-2.9}
\liminf_{n\to\infty}\inte |v_{1,n}|^4dx\geq \delta\quad\hbox{and}\quad\liminf_{n\to\infty}\inte |v_{2,n}|^4dx\geq \delta,
\end{equation}
where $\delta>0$ is independent of $n$.
Combining \eqref{s-2.10} and \eqref{s-2.9} then yields that
\[
\begin{aligned}
e(2,a)&\leq \lim_{n\to\infty}\Big[\lambda F_{2,a}(v_{1,n})+(1-\lambda)F_{2,a}(v_{2,n})\Big]\\
&=\lim_{n\to\infty}\Big[F_{2,a}\big(\sqrt{\lambda}v_{1,n}\big)+ F_{2,a}\big(\sqrt{1-\lambda} v_{2,n}\big)\\
&\qquad\quad\quad-\frac{|a|\lambda(1-\lambda)}{2}\inte (|v_{1,n}|^4+|v_{2,n}|^4)dx\Big]\\
&\leq  e_{\lambda}(2,a)+e_{1-\lambda}(2,a)-|a|\lambda(1-\lambda)\delta\\
&<e_{\lambda}(2,a)+e_{1-\lambda}(2,a),
\end{aligned}
\]
and the claim \eqref{s-2.6} is hence proved.

Based on \eqref{s-2.6}, we now prove that the dichotomy cannot occur. Indeed, taking a subsequence if necessary, we may assume that
\begin{equation}\label{s-2.8}
\lim_{n\to\infty}\inte f_n(x)dx=\alpha>0.
\end{equation}
By contradiction, suppose that the dichotomy occurs. Then there exist $\beta\in (0,\alpha)$, $0<2R_0<R_n\to+\infty$ as $n\to+\infty$ and $\{y_n\}\subset \R^2$ such that for any $\eps>0$, we have
\begin{equation}\label{2.18}
\Big|\int_{B_{R_0}(y_n)} f_n(x)dx-\beta\Big|\leq \eps,\,\ \Big|\int_{\R^2/B_{2R_n}(y_n)}f_n(x)dx-(\alpha-\beta)\Big|\leq \eps\,\ \hbox{as $n\to\infty$}.
\end{equation}
Set $$u_{1n}(x):=\psi\big(\frac{x-y_n}{R_0}\big)u_n(x),\,\ \,\  u_{2n}(x):=\phi\big(\frac{x-y_n}{R_n}\big)u_n(x),$$ where $0\leq\psi(x)\leq 1$ and $0\leq\phi(x)\leq 1$ are smooth cut-off functions and satisfy:  $\psi\equiv 1$ if $|x|\leq 1$, and $\psi\equiv 0$ if $|x|\geq 2$; $\phi\equiv 0$ if $|x|\leq 1$, and $\phi\equiv 1$ if $|x|\geq 2$.
By the definition of $\{u_{1n}\}$ and $\{u_{2n}\}$, we have
\begin{equation}\label{s-2.11}
\hbox{dist(supp $u_{1n}$, supp $u_{2n}$)$\to+\infty$\,\ as\,\  $n\to\infty$.}
\end{equation}
We claim that as $n\to\infty$,
\begin{equation}\label{s-2.12}
\begin{aligned}
&\Big|\|u_{1n}\|^2_{H_A^1}-\beta\Big|\leq C\varepsilon,\quad\Big|\|u_{2n}\|^2_{H_A^1}-(\alpha-\beta)\Big|\leq C\varepsilon,\\
&\|u_n-u_{1n}-u_{2n}\|_{H_A^1}\leq C\varepsilon,
\end{aligned}
\end{equation}
where we denote
\begin{equation}\label{s-2.26}
H_A^1(\R^2,\C):=\big\{u\in L^2(\R^2,\C):\,(\nabla-ix^\perp) u\in L^2(\R^2,\C) \big\}
\end{equation}
associated with the norm
\[
\|u\|_{H^1_A}:=\Big(\inte (|(\nabla-ix^\perp) u|^2+|u|^2)dx\Big)^{\frac 12}.
\]

We now prove the claim (\ref{s-2.12}). Actually, we obtain from \eqref{s-2.8} and \eqref{2.18} that as $n\to\infty$,
\[
\begin{aligned}
&\Big|\|u_{1n}\|^2_{H_A^1}-\beta\Big|\\
\leq &\Big|\|u_{1n}\|^2_{H_A^1}-\int_{B_{R_0}(y_n)} f_n(x)dx\Big|+\Big|\int_{B_{R_0}(y_n)} f_n(x)dx-\beta\Big|\\
\leq &\Big|\int_{B_{2R_0}(y_n)/B_{R_0}(y_n)}\Big\{\big|(\nabla-ix^\perp) \big[\psi(\frac{x-y_n}{R_0})u_{n}\big]\big|^2+|\psi(\frac{x-y_n}{R_0})u_{n}|^2\\
&-\big|(\nabla-ix^\perp) u_{n}\big|^2-|u_{n}|^2\Big\}dx\Big|+\eps\\
\leq& \int_{B_{2R_0}(y_n)/B_{R_0}(y_n)}\Big(3|(\nabla-ix^\perp) u_{n}|^2+2\Big|\nabla \psi\Big(\frac{x-y_n}{R_0}\Big)\Big|^2|u_{n}|^2+2|u_{n}|^2\Big)dx+\eps\\
\leq& \int_{B_{2R_n}(y_n)/B_{R_0}(y_n)}C\big(|(\nabla-ix^\perp) u_{n}|^2+|u_{n}|^2\big)dx+\eps\leq  C\eps.
\end{aligned}
\]
%This further implies that as $n\to\infty$
%\[
%\big|\|u_{1n}\|^2_{H_A^1}-\beta\big|\leq  C\eps.
%\]
Similarly, we can also get that as $n\to\infty$,
\[
\big|\|u_{2n}\|^2_{H_A^1}-(\beta-\alpha)\big|\leq  C\eps.
\]
Following \eqref{s-2.8} and \eqref{2.18}, we also have
\[
\begin{aligned}
&\|u_n-u_{1n}-u_{2n}\|_{H_A^1}^2\\
=&\int_{B_{2R_n}(y_n)/B_{R_0}(y_n)}\Big|\big(\nabla-ix^\perp\big)\Big(1-\psi(\frac{x-y_n}{R_0})-\phi(\frac{x-y_n}{R_n})\Big)u_n
\Big|^2\\
&+\Big|\Big(1-\psi(\frac{x-y_n}{R_0})-\phi(\frac{x-y_n}{R_n})\Big)u_n\Big|^2dx\\
\leq &\int_{B_{2R_n}(y_n)/B_{R_0}(y_n)} C\Big[|(\nabla-ix^\perp)u_n|^2+\Big|\nabla\big[\psi(\frac{x-y_n}{R_0})+\phi(\frac{x-y_n}{R_n})\big]\Big|^2|u_n|^2+|u_n|^2\Big]dx\\
\leq &\int_{B_{2R_n}(y_n)/B_{R_0}(y_n)}Cf_n(x)dx\leq C\eps\ \ \hbox{as}\,\ n\to\infty.
\end{aligned}
\]
We now conclude from above that the claim \eqref{s-2.12} holds true.

Since the embedding $H_A^1(\R^2,\mathbb{C})\hookrightarrow L^p(\R^2,\mathbb{C})$ $(2\leq p<\infty)$ is continuous, cf. \cite[Proposition 2.1]{EL}, we derive from \eqref{s-2.12} that
\begin{equation}\label{s-2.13}
\big|F_{2,a}(u_n)-F_{2,a}(u_{1n})-F_{2,a}(u_{2n})\big|\leq C\varepsilon,
\end{equation}
\begin{equation}\label{s-2.14}
\Big|\inte \big(|u_n|^p-|u_{1n}|^p-|u_{2n}|^p\big)dx\Big|\leq C\varepsilon\quad\hbox{for any $p\in [2,\infty)$}.
\end{equation}
Moreover, following \eqref{GNineq}, \eqref{s-2.19} and \eqref{s-2.12}, there exists a constant $\delta>0$ such that for $i=1,2$,
\begin{equation}\label{s-2.15}
\begin{split}
F_{2,a}(u_{in})&\geq (1-\frac{a}{a^*})\inte |(\nabla-ix^\perp) u_{in}|^2dx\\
&\geq \frac{2(a^*-a)}{3a^*}\|u_{in}\|^2_{H_A^1(\R^2)}\geq \delta\quad \hbox{for all $n$, }
\end{split}
\end{equation}
which then implies that there exist subsequences of $\{u_{1n}\}$ and $\{u_{2n}\}$, still denoted by $\{u_{1n}\}$ and $\{u_{2n}\}$, satisfying $\lim_{n\to \infty}\|u_{1n}\|^2_{L^2(\R^2)}=\theta$ and $\lim_{n\to \infty}\|u_{2n}\|^2_{L^2(\R^2)}=1-\theta$ for some constant $\theta\in(0,1)$ by applying (\ref{s-2.14}).
Since $e_\lambda(2,a)$ is continuous in $\lambda$, we conclude from \eqref{s-2.14} and \eqref{s-2.15} that
\[
\lim_{n\to\infty} F_{2,a}(u_{1n})\geq e_{\theta}(2,a),\, \lim_{n\to\infty} F_{2,a}(u_{2n})\geq e_{1-\theta}(2,a),
\]
which then yield that $e(2,a)\ge e_{\theta}(2,a)+e_{1-\theta}(2,a)$ holds for some constant $\theta\in(0,1)$ in view of (\ref{s-2.13}). This however contradict to
\eqref{s-2.6}. Therefore,  the dichotomy cannot occur.

Since both the vanishing and the dichotomy cannot occur, we now conclude that only the compactness
occurs for the sequence $\{f_n\}$, i.e., there exists $\{x_n\}\subset \R^2$ such that for any fixed $\varepsilon>0$,
\begin{equation}\label{s-2.16}
\int_{\R^2/B_R(x_n)}f_n(x)dx\leq \varepsilon,\quad\hbox{if $R=R(\varepsilon)>0$ is large enough.}
\end{equation}
Since $\{u_n\}$ is bounded in $H_A^1(\R^2,\mathbb{C})$, we obtain from (\ref{s-2.3}) and \eqref{s-2.16} that $\{u_n(\cdot+ x_n)\}$ is relatively compact in $L^2(\R^2)$ by \cite[Proposition 2.1]{EL} and thus relatively compact in $L^p(\R^2)$ $(2<p<+\infty)$  by H\"{o}lder inequality. Note that
$$\tilde u_n(x):=\hat{u}_n(x)e^{-ix_n\cdot x^{\perp}}=u_n(x+ x_n)e^{-ix_n\cdot x^{\perp}}$$
satisfies
\begin{equation}\label{s-2.21}
\inte |(\nabla-ix^\perp)\tilde u_n|^2dx=\inte \big|\big(\nabla -ix^{\perp}\big)\big(\hat u_ne^{-ix_n\cdot x^{\perp}}\big)\big|^2dx=\inte \big|(\nabla-ix^\perp)u_n\big|^2dx.
\end{equation}
%Setting $$, we then have
%\[
%\inte |(\nabla-ix^\perp)u_n|^2dx=\inte |(\nabla-ix^\perp)\tilde u_n|^2dx.
%\]
Moreover, since $|\tilde u_n(x)|=|u_n(x+x_n)|$ holds in $\R^2$, by \eqref{s-2.8} and \eqref{s-2.16} we now derive from above that there exists a subsequence, still denoted by $\{\tilde u_n\}$, of $\{\tilde u_n\}$ satisfying
\begin{equation}\label{s-2.28}
\lim_{n\to\infty}\tilde u_n=u_0\,\ \hbox{strongly in $L^p(\R^2)$ ($2\leq p<\infty$) for some $u_0\in \mathcal{H}$,}
\end{equation}
and
\[
\liminf_{n\to\infty}\inte \big|(\nabla-ix^\perp)u_n\big|^2dx
=\liminf_{n\to\infty}\inte \big|(\nabla-ix^\perp)\tilde u_n\big|^2dx
\geq \inte \big|(\nabla-ix^\perp)u_0\big|^2dx.
\]
Therefore,  $e(2,a)$ admits at least one minimizer $u_0$ in the case $a\in(-a^*,0)$.

(2). We now consider the case $a\in (0,\infty)$.
By selecting a suitable eigenvalue function  of \eqref{s-2.17} as a test function, we first prove that
\begin{equation}\label{2.10M}
e(2,a)\leq 2,\ \ a\in (0,\infty).
\end{equation}
Towards the proof of (\ref{2.10M}), motivated by \cite{A,AB}, we define the following
lattice
\begin{equation}\label{2.10}
\ell:=v(\mathbb{Z}\oplus e^{\frac{2i\pi}{3}}\mathbb{Z})=v(m+ne^{\frac{2i\pi}{3}}),
\end{equation}
where $m,n\in \mathbb{Z}$ are arbitrary integers, and the constant $v\in\R^+$ is chosen such that $\frac{\sqrt{3}v^2}{2}>\pi$. Let $Q$ be a unit cell (i.e., a regular hexagon) centered at the origin, whose area satisfies $|Q|=\frac{\sqrt{3}v^2}{2}$. Here $\frac{1}{|Q|}$ denotes an average spatial density of the point in the lattice.
Choose the following test function
\begin{equation}\label{s-2.22}
\psi_R(z)=A_Re^{-\frac{1}{2}|z|^2}\prod_{j\in\ell\bigcap B_R(0)}(z-j),
\end{equation}
where we identify any complex number $z=x_1+ix_2$ with the point  $x=(x_1,x_2)\in\R^2$, and the Lebesgue measure $dz$ with $dxdy$, respectively. Here $A_R>0$ is chosen properly such that $\|\psi_R\|_{L^2(\R^2)}=1$. Obviously, $\psi_R\in S$, where $S$  is defined by \eqref{s-2.18}, and $(2,\psi_R)$ is the first eigenpair of (\ref{s-2.17}). We hence obtain from  (\ref{s-2.17}) that
\begin{equation}\label{s-2.22M}
\inte \big|(\nabla -ix^{\perp})\psi_R\big|^2dx=2.
\end{equation}
Moreover, it follows from \cite[Theorem 5.2]{A} that as $R\to\infty$,
\begin{equation}\label{s-2.24}
|\psi_R(z)|\to \psi(z)=\frac{1}{\sqrt{\pi}\sigma}|\eta(z)|e^{-\frac{|z|^2}{2\sigma^2}}\quad\hbox{in}\,\ L^p\big(\R^2,(1+|z|^2)dz\big),\ \ p\geq 1,
\end{equation}
where
\begin{equation}\label{s-2.25}
\frac{1}{\sigma^2}=1-\frac{\pi}{|Q|}\to 0\,\ \hbox{as $|Q|\to\pi$},
\end{equation}
and $|\eta(z)|$ is periodic and vanishes at each point of the lattice $\ell$, see more details in \cite[Theorem 5.1]{A}. We should remark that the main feature of the periodic lattice $\ell$ is to
modify the decay of the Gaussian from $e^{-\frac{|z|^2}{2}}$ to $e^{-\frac{|z|^2}{2\sigma^2}}$, where $\sigma>0$ depends on
the area of $Q$ through the relation \eqref{s-2.25}.
%Moreover, $|\eta(z)|$ satisfies
%\[
%-\Delta (ln |\eta|)=2\pi\delta_0-\frac{\pi}{2|Q|}\,\ \hbox{in Q},
%\]
%with periodical boundary condition.
Therefore, we have
\begin{equation}\label{2.29}
\begin{split}
 \quad\lim_{R\to\infty}\inte |\psi_R(z)|^4dz
&= \inte\frac{1}{\pi^2\sigma^4}|\eta(z)|^4e^{-\frac{2|z|^2}{\sigma^2}}dz\\
&=\inte\frac{1}{\pi^2\sigma^2}|\eta(\sigma z)|^4e^{-2|z|^2}dz\\
&\leq \frac{C}{\pi^2\sigma^2}\inte e^{-2|z|^2}dz\leq\frac{C}{\sigma^2}\to 0\,\ \hbox{as} \,\ \sigma\to\infty,
\end{split}
\end{equation}
and thus
\[
\begin{aligned}
e(2,a)&\leq \lim_{\sigma\to\infty,\,R\to\infty}F_{2,a}(\psi_R)\\
&=2+\lim_{\sigma\to\infty,\,R\to\infty}\inte |\psi_R|^4dz
\leq 2, \ \ \mbox{if}\ \ a\in (0,\infty),
\end{aligned}
\]
i.e.,  the estimate (\ref{2.10M}) holds true.

We now assume that $e(2,a)$ admits a minimizer  $u_a$ for any $a\in (0,\infty)$. We then obtain from \eqref{s-2.19} that
\[
 \inte \big|(\nabla-ix^\perp)u_a\big|^2dx\geq 2 \inte |u_a|^2dx=2,
\]
which and (\ref{2:newF}) thus yield that
\[
e(2,a)> 2, \ \ \mbox{if}\ \ a\in (0,\infty),
\]
a contradiction to (\ref{2.10M}). We therefore conclude that there is no minimizer for $e(2,a)$ in the case $a\in (0,\infty)$. This therefore completes the proof of Theorem \ref{thm1.1} \qed

%(2.B). If $a=-a^*$, one can get from (\ref{s-2.1}) that $F_{2,-a^*}(u)\geq0$ for any $u\in\mathcal{H}$ satisfying $\|u\|_{L^2(\R^2)}^2=1$. However, one can also get that $e(2,-a^*)\le 0$ by taking $u_\tau(x)=\frac{1}{\sqrt{a^*}}\tau w(\tau x)$ as a test function and then letting $\tau\to\infty$. We hence have $e(2,-a^*)= 0$.
%
%
%If $e(2,-a^*)$ admits a minimizer $u_a$, then we obtain from (\ref{GNineq}) and (\ref{Diam}) that
%\begin{eqnarray*}
%0&=&\int_{\mathbb{R}^2} |(\nabla-ix^{\perp})u_a|^2dx-\frac{a^*}{2}\inte|u_a|^4dx\\
%&\geq& \inte \big|\nabla |u_a|\big|^2dx-\frac{a^*}{2}\inte|u_a|^4dx\ge 0.
%\end{eqnarray*}
%This further implies that
%\[\int_{\mathbb{R}^2} |(\nabla-ix^{\perp})u_a|^2dx=\inte \big|\nabla |u_a|\big|^2dx=\frac{a^*}{2}\inte|u_a|^4dx,
%\]
%which is impossible. We therefore have the nonexistence of minimizers for $e(2,-a^*)$.
%
%
%(2.C). Finally, if $a\in(-\infty,-a^*)$, then one can have $e(2,a)\le \lim_{\tau\to\infty}F_{2,a}\big(\tau w(\tau x)\big)=-\infty$, which implies the nonexistence of minimizers in this case.

\vskip 0.2truein
\noindent{\bf Proof of Theorem \ref{thm1.1*}.} Without loss of generality, as before we may assume $A=1$ so that $\Omega^*=2$. If $a\in (-\infty,-a^*]$, the argument of proving \cite[Theorem 2.3]{GLY} then gives the nonexistence of minimizers for $e(2,a)=e(\Omega^*,a)$. Therefore, we only   consider the case where $a\in (-a^*, +\infty)$.
	
(\rm i). Assume $a\in (-a^*, +\infty)$ and $W(x)>B$ in $\R^2$. For simplicity we may assume $B=0$ so that $W(x)>0$ in $\R^2$.

We first consider the case $a\in(-a^*,0]$. In this case, we claim that
\begin{equation}\label{s-2.23}
e(2,a)=e^\star(2,a),\quad\hbox{where $e^\star(2,a)$ denotes the value of $e(2,a)$ at $W\equiv 0$.}
\end{equation}
Recall from Theorem \ref{thm1.1} (i) that  $e^\star(2,a)$ admits a minimizer $\hat u(x)\in \mathcal{H}$. Under the assumptions of $W(x)$, we have
\begin{equation}\label{2.12}
e(2,a)\geq e^\star(2,a).
\end{equation}
On the other hand, taking the following test function
\[
\hat u_k(x):=\hat u(x+x_k)e^{-ix_k^\perp\cdot x},\quad\hbox{where}\,\ |x_k|\to\infty\ \ \mbox{as}\ \ k\to\infty,
\]
we have
\[
e(2,a)\leq F_{2,a}(\hat u_k)\to e^\star(2,a)\ \ \mbox{as}\ \ k\to\infty.
\]
We then obtain from \eqref{2.12} and above that the claim \eqref{s-2.23} is true. Following this estimate, since $W(x)>0$ in $\R^2$, as in \cite[Theorem 1]{GS} one can further derive the nonexistence of minimizers for $e(2,a)$ in the case $a\in(-a^*,0]$.

We next consider the case $a\in(0,+\infty)$. Similar to Theorem \ref{thm1.1} (ii), we take a test function as follows:
\begin{equation}\label{new3}
\psi_{R,\sigma}(z):=\psi_R(z+z_\sigma)e^{-ix_\sigma^\perp\cdot x},
%=A_Re^{-\frac{1}{2}|z+z_\sigma|^2+ix_\sigma^\perp\cdot x}\prod_{j\in\ell\bigcap B_R(0)}(z+z_\sigma-j),
\end{equation}
where $\psi_R(z)$ and the parameter $\sigma$ are defined by \eqref{s-2.22} and \eqref{s-2.25}, respectively. Here $x_\sigma=(x_{1,\sigma},x_{2,\sigma})$ and $z_\sigma=x_{1,\sigma}+ix_{2,\sigma}$ satisfy $|x_\sigma|=|z_\sigma|=\sigma^2$.
One can check that $\psi_{R,\sigma}(z)\in S$ satisfies
\[
\inte \big|(\nabla-ix^{\perp})\psi_{R,\sigma}\big|^2dx=\inte \big|(\nabla-ix^{\perp})\psi_{R}\big|^2dx=2,
\]
and
\[
\inte |\psi_{R,\sigma}|^2dx=\inte |\psi_R|^2dx=1.
\]
Note from \eqref{2.29} that
\[
\lim_{R\to\infty}\inte |\psi_{R,\sigma}|^4dx=\lim_{R\to\infty}\inte |\psi_R|^4dx\to 0\quad\hbox{as}\,\  \sigma\to+\infty.
\]
%where $\frac{1}{\sigma^2}=1-\frac{\pi}{|Q|}$ and $|Q|>\pi$ denotes the area of $Q\in \ell$ for $\ell$ defined by \eqref{2.10}.
Moreover, because $|z_\sigma|=\sigma^2$, we derive from \eqref{s-2.24} that
\[
\begin{aligned}
&\lim_{R\to\infty}\inte W(z)|\psi_R(z)|^2dz
=\inte \frac{1}{\pi \sigma^2}W(z+z_\sigma)|\eta(z)|^2e^{-\frac{|z|^2}{2\sigma^2}}dz\\
=&\int_{|z|<\sigma^{\frac{1}{2}}}\frac{1}{\pi }W(\sigma z+z_\sigma)|\eta(\sigma z)|^2e^{-\frac{|z|^2}{2}}dz+\int_{|z|\geq \sigma^{\frac{1}{2}}}\frac{1}{\pi }W(\sigma z+z_\sigma)|\eta(\sigma z)|^2e^{-\frac{|z|^2}{2}}dz\\
:=&I+II.
\end{aligned}
\]
Since $\lim_{|x|\to\infty}W(x)=0$ and $|\eta(z)|\in L^\infty(\R^2)$, we get that $|I|\leq C(\sigma)\to 0$ as $\sigma\to\infty$. Moreover, the exponential decay of $e^{-\frac{|z|^2}{2}}$ yields that $|II|\leq C(\sigma)\to 0$ as $\sigma\to\infty$. Therefore, we obtain from above that
\[
\begin{aligned}
&\lim_{R\to\infty}\inte W(x)|\psi_{R,\sigma}(x)|^2dx=\lim_{R\to\infty}\inte W(x+x_\sigma)|\psi_R(x)|^2dx\to 0\,\ \hbox{as $\sigma\to\infty$}.
\end{aligned}
\]
We thus conclude from \eqref{2.12} and above that $e(2,a)=e^\star(2,a)$. Since $W(x)>0$ in $\R^2$, as before one can further  obtain  the nonexistence of minimizers for $e(2,a)$ in the case $a\in(0,+\infty)$. This completes the proof of Theorem \ref{thm1.1} (\rm i).

\vskip 0.1truein
(\rm ii). Let $a\in (-a^*, +\infty)$ and $W(x)<B$ in $\R^2$. For simplicity we may assume $B=0$  so that $W(x)<0$ in $\R^2$. Following the definition of $a_\ast$
in \eqref{ast}, we can show that $a_\ast>0$ by using \cite[Theorem 2]{LP}. Indeed, for $A=1$ and $B=0$, we can take $u=e^{-\frac{|x|^2}{2}}$ as a test function to obtain that
\[
e(2,a)\leq F_{2,a}(u)=2+\inte \frac{W(x)e^{-|x|^2}}{\pi}dx+\frac{a}{2}\inte \frac{e^{-2|x|^2}}{\pi^2}dx.
\]
Since $W(x)<0$ in $\R^2$, we get that $e(2,a)<2$ if $a>0$ is small enough and thus $a_\ast>0$ holds. 

(a). We first consider $a\in(-a^*,a_\ast)$. By the definition of $a_\ast$, it is obviously that
\begin{equation}\label{s-2.29}
e(2,a)<2\quad\hbox{for any $a\in (-a^*, a_\ast)$.}
\end{equation}
Similar to Theorem \ref{thm1.1} (\rm ii) we shall prove the existence of minimizers for $e(2,a)$ by the concentration-compactness lemma.
Towards this aim, we denote $\{u_n\}$ a minimizing sequence of $e(2,a)$, i.e., $\|u_n\|^2_2=1$ and $e(2,a)=\lim_{n\to\infty}F_{2,a}(u_n)$.
As in Theorem \ref{thm1.1} (ii), we denote
\begin{equation}\label{s-2.3N}
f_n(x):=|(\nabla -ix^{\perp})u_n|^2+|u_n|^2\ge 0\ \ \mbox{in}\, \ \R^2.
\end{equation}
One can check that the nonnegative sequence $\{f_n\}$ is bounded uniformly in $L^1(\R^2)$.

Following \eqref{s-2.29}, the same argument of Theorem \ref{thm1.1} (ii) yields that the vanishing cannot occur for $\{f_n(x)\}$ defined in (\ref{s-2.3N}).
Moreover, applying \eqref{s-2.29}, one can also derive that \eqref{s-2.6} holds for any $a\in(-a^*,+\infty)$. Following this estimate, the same argument of Theorem \ref{thm1.1} (\rm ii) yields that the dichotomy cannot occur  for $\{f_n(x)\}$. By the concentration-compactness lemma, we therefore conclude that the compactness holds for $\{f_n(x)\}$.

Finally, the same argument of proving \eqref{s-2.28} yields that $\tilde u_n=u_n(\cdot+x_n)e^{ix_n\cdot x^\perp}$ converges to some $u_0\in \mathcal{H}$ in $L^p(\R^2)$ $(2\leq p<\infty)$ as $n\to\infty$. To obtain a true minimizer of $e(2,a)$, we next show that $\{x_n\}$ is bounded uniformly in $\R^2$. On the contrary, suppose  $$|x_n|\to+\infty\,\ \hbox{as}\,\ n\to\infty.$$
Then up  to a subsequence if necessary, we have
\[
\lim_{n\to\infty}\inte W(x)|u_n(x)|^2dx=\lim_{n\to\infty}\inte W(x+x_n)|u_n(x+x_n)|^2dx=0.
\]
By the definition of  $e^*(2,a)$ defined in \eqref{s-2.23}, we thus have
\[
e(2,a)=\lim_{n\to\infty}F(2,a)(u_n)
\geq \inte\big(|(\nabla-ix^\perp)u_0|^2+\frac{a}{2}|u_0|^4\big)dx\geq e^\star(2,a),
\]
which however contradicts to \eqref{s-2.23}. Once $\{x_n\}$ is bounded uniformly in $\R^2$, we may assume that  up  to a subsequence if necessary, $x_n\to y_0$ for some $y_0\in\R^2$ as $n\to\infty$. Moreover, since $W(x)\in C(\R^2)$ and $\lim_{|x|\to \infty}W(x)=0$, we obtain that $W(x+x_n)$ converges to $W(x+y_0)$ uniformly in $\R^2$ as $n\to\infty$. Therefore, we get that
\[
\begin{aligned}
\lim_{n\to\infty}\inte W(x)|u_n(x)|^2dx=&\lim_{n\to\infty}\inte W(x+x_n)|u_n(x+x_n)|^2dx\\
=&\lim_{n\to\infty}\inte W(x+y_0)|u_n(x+x_n)|^2dx\\
=&\inte W(x+y_0)|u_0(x)|^2dx,
\end{aligned}
\]
where we used the fact that $W(x)\in L^\infty(\R^2)$ and $ \lim_{n\to\infty}\inte|u_n(\cdot+x_n)|^2dx=\inte |u_0|^2dx$ in the last equality.
This further implies that
\[
\begin{aligned}
e(2,a)=\lim_{n\to\infty}F_{2,a}(u_n)
&\geq  \inte\big(|(\nabla  -ix^{\perp})u_0|^2+W(x+y_0)|u_0|^2+\frac{a}{2}|u_0|^4\big)dx\\
&=F_{2,a}\big(u_0(\cdot-y_0)e^{iy_0\cdot x^\perp}\big),
\end{aligned}
\]
i.e., $u_0(\cdot-y_0)e^{iy_0\cdot x^\perp}$ is a minimizer of $e(2,a)$ for the case $a\in(-a^*,a_\ast)$, and we are done.

(b). We next consider $a\in(a_\ast,+\infty)$. By the definition of $a_\ast$, we obtain that $e(2,a)\geq 2$ for any $a\in(a_\ast,+\infty)$. On the other hand, taking the same test function as in \eqref{new3}, we can obtain that $e(2,a)\leq 2$ and hence $e(2,a)\equiv2$ for any $a\in(a_\ast,+\infty)$. Therefore, we can prove that there is no minimizers for $e(2,a)=2$ in this case. Indeed, assume that there is some $u_0\in H^1(\R^2,\mathbb{C})$ satisfying $\|u_0\|^2_{L^2(\R^2)}=1$ such that $e(2,a)=2$. We then have $2=e(2,a)\geq e(2,a_\ast)+(a-a_\ast)\inte |u_0|^4dx=2$, and thus $u_0\equiv 0$, which contradicts with $\|u_0\|^2_{L^2(\R^2)}=1$. This completes the proof of Theorem \ref{thm1.1*}.  \qed

\vskip 0.2truein
\noindent{\bf Proof of Corollary \ref{cor-1.3}.}
Without loss of generality, as before we may assume $A=1$ so that $\Omega^*=2$, and $B=0$.
Following Theorem \ref{thm1.1*} (ii), we only need to prove that $a_\ast=+\infty$, or equivalently, that $e(2,a)<2$ for any $a>0$.
Taking $\psi_R$ as a test function, where  $\psi_R$ is defined by (\ref{s-2.22}), we  have
\[
\begin{aligned}
e(2,a)&\leq \lim_{R\to\infty}F_{2,a}(\psi_R)=2+\lim_{R\to\infty}\inte \Big(W(z)|\psi_R(z)|^2+\frac{a}{2} |\psi_R(z)|^4\Big)dz\\
&=2+\inte \Big(W(z)\frac{|\eta(z)|^2e^{-\frac{|z|^2}{\sigma^2}}}{\pi\sigma^2}+\frac{a}{2} \frac{|\eta(z)|^4e^{-\frac{2|z|^2}{\sigma^2}}}{\pi^2\sigma^4}\Big)dz\\
&=2+\inte \Big(W(\sigma z)\frac{|\eta(\sigma z)|^2e^{-|z|^2}}{\pi}+\frac{a}{2} \frac{|\eta(\sigma z)|^4e^{-2|z|^2}}{\pi^2\sigma^2}\Big)dz,
\end{aligned}
\]
where $\eta(z)$ is as in \eqref{s-2.24}.
Since $|\eta(z)|$ is a periodic function, we obtain from \cite{Al} that $|\eta(\sigma z)|^2$ and $|\eta(\sigma z)|^4$ converge $L^\infty-$weak star to the constants $\bbint|\eta(z)|^2dz>0$ and $\bbint|\eta(z)|^4dz>0$ as $\sigma\to\infty$,
respectively, where $\bbint |\eta(z)|^pdz=\frac{1}{|Q|}\int_Q |\eta( z)|^pdz$ denotes the average of $|\eta(z)|^p$ in $Q$.
Following \eqref{new2}, $\limsup_{|x|\to\infty }W(x)|x|^s<0$ holds for $0<s<2$, we obtain that there exists a small $\delta>0$ such that
\[
\begin{aligned}
&\quad\inte \Big(W(\sigma z)\frac{|\eta(\sigma z)|^2e^{-|z|^2}}{\pi}+\frac{a}{2} \frac{|\eta(\sigma z)|^4e^{-2|z|^2}}{\pi^2\sigma^2}\Big)dz\\
&\leq \int_{\R^2/B_1(0)}\frac{-\delta |\eta(\sigma z)|^2e^{-|z|^2}}{\pi\sigma^s|z|^s}dz+\frac{C}{\sigma^2}\\
&\leq  -\frac{C(\delta)}{\sigma^s}+\frac{C}{\sigma^2}<0\quad\hbox{as}\,\ \sigma\to\infty,
\end{aligned}
\]
which then gives that
\[
e(2,a)<2 \quad\hbox{for any}\,\ a>0,
\]
and hence $a_\ast=+\infty$ and the proof is completed in view of Theorem \ref{thm1.1*} (\rm ii).  \qed

%We first claim that Following Theorem \ref{thm1.1} (i), $e^*(2,a)$ admits at least one minimizer denoted by $\hat u(x)$. Since $W(x)<0$ in $\R^2$, we then have
%\[
%e(2,a)\leq F_{2,a}(\hat u)=e^*(2,a)+\inte W(x)|\hat u(x)|^2dx<e^*(2,a),
%\]
%which implies that (\ref{s-2.23}) holds for $a\in(-a^*,0]$.
%However, if $a\in(0,+\infty)$, then we use the same test function $\psi_R$ defined in \eqref{s-2.22}.

\section{Limit behavior of minimizers as $\Omega \searrow 0$}

For any fixed $a\in[0,+\infty)$ or a.e. $a\in (-a^*,0)$, in this section we analyze the refined limit behavior of minimizers for $e(\Omega,a)$ as $\Omega>0$ approaches to zero, where the trap $V(x)\ge 0$ satisfies all assumptions of Theorem \ref{thm1.2}.

Towards the above purpose, we denote $u_n$ a complex-valued minimizer of $e(\Omega_n,a)$, where $\Omega_n\searrow 0$ as $n\to\infty$, so that $u_n$  satisfies the following Euler-Lagrange equation
\begin{equation}\label{3.2}
-\Delta u_n+\big(V(x)-\mu_n\big)u_n+i\, \Omega_n \, (x^{\perp}\cdot \nabla u_n)+a|u_n|^2u_n=0\ \ \mbox{in}
\, \  \R^2,
\end{equation}
where $\mu_n\in\R$ is a Lagrange multiplier and satisfies
\begin{equation}\label{3.4}
\mu_n=e(\Omega_n,a)+\frac{a}{2}\inte|u_n|^4dx.
\end{equation}
Under the assumptions of Theorem \ref{thm1.2}, recall from Theorem B that $u_0>0$ is the unique positive minimizer of $e(a)$. Denote $\mu_0\in\R$ the Lagrange multiplier  associated to $u_0$, so that $(u_0,\mu _0)$ satisfies
\begin{equation}\label{3.4MM}
-\Delta u_0+V(x)u_0=\mu_0u_0-au_0^3\ \ \mbox{in}
\, \  \R^2.
\end{equation}
We start with the following estimates of $e(\Omega_n,a)$ and  $\mu_n$ as $n\to\infty$.

\begin{lem}\label{lem-2.1}
Under the assumptions of Theorem \ref{thm1.2}, let $u_n$ be a complex-valued minimizer of $e(\Omega_n,a)$, where $a\in (-a^*,+\infty)$ is fixed, and $\Omega_n\searrow 0$ as $n\to\infty$. Then we have
\begin{enumerate}
\item The energy $e(\Omega_n,a)$ satisfies
\begin{equation}\label{2.8}
\lim_{n\to\infty}e(\Omega_n,a)=e(a),
\end{equation}
where the energy $e(a)$ is defined in (\ref{2.2}).

\item There exists $\theta_n\in [0,2\pi)$ such that $u_n$ satisfies
\begin{equation}\label{2.3}
\lim_{n\to\infty}u_ne^{i\theta_n}=u_0\,\ \hbox{strongly in}\,\ H^1(\R^2,\mathbb{C}),
\end{equation}
where $u_0>0$ is the unique positive minimizer of $e(a)$, and the Lagrange multiplier $\mu_n=\mu(\Omega_n,a)$ of (\ref{3.2}) satisfies
\begin{equation}\label{2.9}
\lim_{n\to\infty}\mu_n=\mu_0,
\end{equation}
where $\mu_0\in\R $ is the Lagrange multiplier of (\ref{3.4MM}).

\item $u_n$ satisfies the following exponential decay
\begin{equation}\label{2.25}
|u_n(x)|,\,|\nabla u_n(x)|\leq C(a)e^{-2|x|}\,\ \hbox{in}\,\ \R^2,
\end{equation}
where the constant $C(a)>0$ depends only on $a$.
\end{enumerate}
\end{lem}

\noindent{\bf Proof.} 1. Taking $u_0$ as a test function, we obtain the following upper bound estimate:
\begin{equation}\label{2.21}
\lim_{n\to\infty}e(\Omega_n,a)\leq \lim_{n\to\infty}F_{\Omega_n,a}(u_0)= e(a).
\end{equation}
On the other hand, we claim that there exists a constant $C>0$, independent of $n$, such that
\begin{equation}\label{2.23}
\inte \big(|\nabla u_n|^2+V(x)|u_n|^2\big)dx\leq C\,\ \hbox{uniformly in}\,\  n.
\end{equation}
Indeed, following \eqref{2.21}, we obtain from \eqref{GNineq} and \eqref{Diam} that
\[
\begin{aligned}
1+e(a)&\geq e(\Omega_n,a)\\
&=\inte\Big(|\nabla u_n|^2-\Omega_nx^\perp\cdot(iu_n,\nabla u_n)+V(x)|u_n|^2+\frac{a}{2}|u_n|^4\Big)dx\\
&\geq\inte\Big(\big|\nabla |u_n|\big|^2+\big[V(x)-\frac{\Omega_n^2}{4}|x|^2\big]|u_n|^2+\frac{a}{2}|u_n|^4\Big)dx\\
&\geq \inte\Big(\big|\nabla |u_n|\big|^2+\frac{a}{2}|u_n|^4\Big)dx-C\\
&\geq\Big(\frac{a^*+a}{2}\Big)\inte |u_n|^4dx-C\quad\hbox{uniformly in $n$.}
\end{aligned}
\]
Since $a\in(-a^*,\infty)$ is fixed, we deduce from above that
\begin{equation}\label{3.1}
\inte |u_n|^4dx\leq C(a)\quad\hbox{uniformly in $n$},
\end{equation}
which further implies that
\[
\begin{aligned}
C(a)&\geq \inte\Big(|\nabla u_n|^2-\Omega_nx^\perp\cdot(iu,\nabla u)+V(x)|u_n|^2\Big)dx\\
&\geq \inte\Big(\frac{1}{2}|\nabla u_n|^2+\big[V(x)-\frac{\Omega_n^2}{2}|x|^2\big]|u_n|^2\Big)dx \\
&\ge \frac{1}{2}\inte\Big[|\nabla u_n|^2+\big(V(x)-C\big)|u_n|^2\Big]dx  \quad\hbox{uniformly in $n$.}
\end{aligned}
\]
This proves the claim \eqref{2.23}.

By the assumptions of $V(x)$, we now deduce from \eqref{2.23} that
\begin{eqnarray*}
\frac{\Omega_n^2}{4}\inte |x|^2|u_n|^2dx &\leq& \frac{\Omega_n^2}{4}\Big[C\inte V(x)|u_n|^2dx+C\Big]\\
&\leq& \frac{\Omega_n^2}{4}C(a)\to 0\quad\hbox{as $n\to\infty$,}
\end{eqnarray*}
which then gives that
\begin{equation}\label{2.22}
\begin{split}
&\lim_{n\to\infty}e(\Omega_n,a)=\lim_{n\to\infty}F_{\Omega_n,a}(u_n)\\
=&\lim_{n\to\infty}\inte \Big(\Big|\nabla u_n-\frac{i\Omega_n x^\perp}{2}u_n\Big|^2+\Big[V(x)-\frac{\Omega_n^2}{4}|x|^2\Big]|u_n|^2+\frac{a}{2}|u_n|^4\Big)dx\\
\geq& \lim_{n\to\infty}\inte \Big(\big|\nabla |u_n|\big|^2+\Big[V(x)-\frac{\Omega_n^2}{4}|x|^2\Big]|u_n|^2+\frac{a}{2}|u_n|^4\Big)dx\\
=&\lim_{n\to\infty}\Big(E_a(|u_n|)-\frac{\Omega_n^2}{4}\inte |x|^2|u_n|^2dx\Big)\geq e(a).
\end{split}
\end{equation}
Therefore, \eqref{2.8} is proved in view of \eqref{2.21} and \eqref{2.22}.

2. For the unique positive minimizer $u_0>0$ of \eqref{2.2}, we choose $\theta_n\in[0,2\pi)$ properly such that
\begin{equation}\label{2.4}
\big\|u_ne^{i\theta_n}-u_0\big\|_{L^2(\R^2)}=\min_{\theta\in [0,2\pi)}\big\|u_ne^{i\theta}-u_0\big\|_{L^2(\R^2)},
\end{equation}
which further implies the following orthogonal condition
\begin{equation}\label{2.4*}
\inte u_0Im (u_ne^{i\theta_n})dx=0\,\ \hbox{for all $n\in \mathbb{N}^+$,}
\end{equation}
where $Im (u)$ denotes the imaginary part of $u$. Taking a subsequence if necessary,  we deduce from \eqref{2.23} and Lemma \ref{2:lem1} that there exists $\hat u\in H^1(\R^2,\mathbb{C})$ such that $u_ne^{i\theta_n}\rightharpoonup \hat u$ weakly in
$H^1(\R^2,\mathbb{C})$ and $u_ne^{i\theta_n}\to\hat u$ strongly in $L^p(\R^2)$ ($2\leq p<\infty$) as $n\to\infty$.
Thus, we get $\|\hat u\|_{L^2(\R^2)}=1$.
Moreover, by the weak lower semi-continuity and Fatou's Lemma, we have
\[
\begin{aligned}
e(a)&=\lim_{n\to\infty}F_{\Omega_n}(u_n)\geq \inte \big(|\nabla \hat u|^2+V(x)|\hat u|^2+\frac{a}{2}|\hat u|^4\big)dx\geq e(a),
\end{aligned}
\]
which therefore implies that $e(a)=E(\hat u)$ and $\hat u=u_0$ in view of \eqref{2.4*}.

Note from  Theorem B that $u_0>0$ is the unique positive minimizer of $e(a)$. We then conclude that the above convergence $u_ne^{i\theta_n}\to u_0$ strongly in $L^p(\R^2)$ ($2\leq p<\infty$) as $n\to\infty$ holds for the whole sequence.  Moreover, since $\lim_{n\to\infty}F_{\Omega_n}(u_n)=e(a)$ and $u_n\rightharpoonup u_0$  weakly in $H^1(\R^2)$, we obtain that $\lim_{n\to\infty}\inte|\nabla u_n|^2dx=\lim_{n\to\infty}\inte|\nabla u_0|^2dx$, and thus \eqref{2.3} holds true.
Finally, one can deduce from \eqref{3.4}--\eqref{2.3} that \eqref{2.9} also holds true.

3.
%By the definition of $\mu_n$, we derive from \eqref{2.8} and \eqref{2.3} that
%\begin{equation}\label{3.7}
%\lim_{n\to \infty} \mu_n=\mu_0.
%\end{equation}
Define $U_n(x) =|u_n(x)|^2\ge 0$ in $\R^2$. We then derive from \eqref{3.2} that
\begin{equation}\label{NEW9}\arraycolsep=1.5pt\begin{array}{lll}
&&
-\displaystyle\frac{1}{2} \Delta U_n+|\nabla u_n|^2-\Omega_n\, x^\bot\cdot(iu_n,\nabla u_n)+\displaystyle\Big[V(x)-\mu _n+a U_n\Big] U_n=0
\quad \mbox{in}\,\ \R^2.\end{array}
\end{equation}
By Cauchy's inequality,
$$|\nabla u_n|^2-\Omega_n\, x^\bot\cdot(iu_n,\nabla u_n)+\frac{\Omega_n^2}{4}|x|^2 U_n \geq 0\ \ \mbox{in} \,\ \R^2,$$
we have
\begin{equation}\label{NEW10}
-\frac{1}{2} \Delta U_n+\Big[V(x)-\frac{\Omega_n^2|x|^2}{4}-\mu _n\Big]U_n+aU_n^2\leq 0 \quad \mbox{in}\,\ \R^2.
\end{equation}
Note that $\lim_{|x|\to \infty} \Big(V(x)-\frac{\Omega_n^2|x|^2}{4}\Big)=+\infty$ and $\lim_{n\to \infty}\mu_n=\mu_0$. By De Giorgi-Nash-Moser theory \cite[Theorem 4.1]{HL},  we derive from above and \eqref{2.3} that
\begin{equation}\label{new6}
0\leq U_n(x)\leq C(a)\quad\hbox{in $\R^2$ uniformly in $n$. }
\end{equation}
It follows from \eqref{new6} that there exists a sufficiently large $R=R(a)>0$ such that $V(x)-\frac{\Omega_n^2|x|^2}{4}-\mu_n-|a| U_n \geq 18$ holds in $\R^2/B_R(0)$ as $n\to\infty$. This then implies from \eqref{NEW10} that
\[
-\frac{1}{2} \Delta U_n+18U_n\leq 0 \quad \mbox{in $\R^2/B_R(0)$ uniformly in $n$}.
\]
By the comparison principle, we obtain from above and \eqref{new6} that
\begin{equation}\label{3.5}
|u_n(x)|\leq C(a)e^{-3|x|}\quad \mbox{in $\R^2$ uniformly in $n$.}
\end{equation}
Based on \eqref{3.5},  the same argument of proving \cite[Lemma 4.2] {GLY} then yields that the gradient estimate
of \eqref{2.25} also holds true, and we are done. \qed

\vskip 0.05truein

Applying \eqref{2.3} and \eqref{2.25}, by a standard bootstrap argument one can further verify from \eqref{3.2} that
\begin{equation}\label{2.24}
\lim_{n\to\infty}u_ne^{i\theta_n}=u_0\,\ \hbox{uniformly in}\,\ L^\infty(\R^2,\mathbb{C}),
\end{equation}
where $u_0>0$ is the unique positive minimizer of $e(a)$, see \cite[Propostion 3.3]{GLY} for the related argument.
We now rewrite
\begin{equation}\label{4.10}
u_ne^{i\theta_n}:=q_n+ir_n,
\end{equation}
where $\theta_n\in [0,2\pi)$ is as in (\ref{2.3}), while $q_n$ and $r_n$ denote the real and imaginary parts of $u_ne^{i\theta_n}$, respectively.
It then follows from (\ref{2.24}) that
\begin{equation}\label{2.28M}
q_n\to u_0 \,\ \hbox{and}\,\  r_n\to 0 \,\ \hbox{uniformly in}\,\ L^\infty(\R^2) \,\ \hbox{as}\,\ n\to\infty.
\end{equation}
Note from \eqref{3.2} that $(q_n,r_n)$ satisfies the following system:
\begin{equation}\label{rev-4:3}
\left\{
\begin{aligned}
\mathcal{L}_n q_n&=\Omega_n (x^\bot\cdot\nabla r_n)
\, \ &\mbox{in}\,\  \R^2,\\
\mathcal{L}_n r_n&=-\Omega_n(x^\bot\cdot\nabla q_n)
\ \ &\mbox{in}\,\  \R^2,\\
\end{aligned}
\right.
\end{equation}
where the operator $\mathcal{L}_n$ is defined by
\begin{equation}\label{2.28}
\mathcal{L}_n:=-\Delta+V(x)-\mu_n+a|u_n|^2.
\end{equation}
In order to obtain some refined estimates of $(q_n,r_n)$, we need to analyze the associated linearized operators $\mathcal{L}$ and $\mathcal{N}$: $M(\R^2,\R)\longmapsto L^2(\R^2,\R)$ defined by
\begin{equation}\label{2.5}
\mathcal{L}:=-\Delta+V(x)-\mu_0+au_0^2,
\end{equation}
\begin{equation}\label{2.6}
\mathcal{N}:=-\Delta+V(x)-\mu_0+3au_0^2,
\end{equation}
where $M(\R^2,\R)$ is denoted by
\[
M(\R^2,\R):=H^2(\R^2,\R)\cap \Big\{u\in L^2(\R^2,\R):\inte V(x)u^2dx< \infty \Big\}.
\]
Since $V(x)\ge 0$ is radially symmetric and satisfies the assumption ($V$), we next address the following analytical properties of  $\mathcal{L}$ and $\mathcal{N}$.

\begin{lem}\label{lem-2}
Under the assumptions of Theorem \ref{thm1.2}, we have
\begin{enumerate}
\item The linearized operator $\mathcal{L}$ satisfies
\begin{equation}\label{2.7}
ker \mathcal{L}=\{u_0\}\,\ \hbox{and}\,\ \langle\mathcal{L}u,u\rangle\geq 0\,\ \hbox{for all}\,\ u\in L^2(\R^2,\R),
\end{equation}
where $u_0>0$ is the unique positive minimizer of $e(a)$.
Moreover, there exists $\rho>0$
such that for all $u\in H^1(\R^2,\R)$,
\begin{equation}\label{1.76}
\langle\mathcal{L}u,u\rangle\geq \rho \|u\|^2_{H^1(\R^2)},\,\  \hbox{if}\,\ \inte u(x)u_0(x)dx=0.
\end{equation}

\item  If $a\neq 0$, then the linearized operator $\mathcal{N}$ is non-degenerate, in the sense that
\begin{equation}\label{2.19}
ker \mathcal{N}=\{0\},
\end{equation}
and $\mathcal{N}^{-1}$: $L^2(\R^2)\longmapsto L^2(\R^2)$ exists and is a continuous linear operator. Moreover, there exists $\hat\rho>0$ such that for all $u\in H^2(\R^2,\R)$,
\begin{equation}\label{2.20}
\|\mathcal{N}u\|_{L^2(\R^2)}\geq \hat\rho \|u\|_{H^2(\R^2)}.
\end{equation}

\item If $a\neq 0$ and $\phi(|x|)\in L^2(\R^2)$ is radially symmetric, then $\psi(x)=\mathcal{N}^{-1}\phi(|x|)\in M(\R^2,\R)$ is also radially symmetric.
\end{enumerate} 	
\end{lem}
{\noindent \bf Proof.} (1). Define the first eigenvalue $\lambda_1$ of $\mathcal{L}$ by
\begin{equation}\label{3.6}
\lambda_1:=\inf\Big\{\inte |\nabla u|^2+(V(x)-\mu_0+a|u_0|^2)|u|^2dx:\,\ u\in M(\R^2,\R),\,\ \|u\|_{L^2(\R^2)}=1\Big\},
\end{equation}
and let $(\lambda_1, \phi)$ be the first eigenpair of $\mathcal{L}$. One can derive easily that $\phi$ does not change sign in $\R^2$.
Since $(0, u_0)$, where $u_0>0$, is an  eigenpair of $\mathcal{L}$, we claim that $\lambda_1=0$. Otherwise, if $\lambda_1<0$, then the eigenfunction $u_0>0$ must be orthogonal to $\phi$, which is impossible in view of the fact that $\phi>0$. Therefore, we conclude from Theorem B that $u_0>0$ is the unique eigenfunction associated to the first eigenvalue $\lambda_1=0$.

Following \cite[Theorem XIII.47 and XIII.67]{RS}, we deduce from above that $\mathcal{L}$ has a non-degenerate ground state and purely discrete spectrum $\sigma(\mathcal{L})=\sigma_d(\mathcal{L})$,  which further imply that
\[
ker \mathcal{L}=\{u_0\}\,\ \hbox{and}\,\ 0\in \sigma_d(\mathcal{L}).
\]
Therefore, for any $u\in H^1(\R^2)$ satisfying $\inte u(x)u_0(x)dx=0$, we have
\[
\langle\mathcal{L}u,u\rangle\geq \rho \|u\|^2_{L^2(\R^2)}.
\]
On the other hand, since $u_0\in L^\infty(\R^2)$, we also reduce from (\ref{3.6}) that for any $u\in H^1(\R^2)$ satisfying $\inte u(x)u_0(x)dx=0$,
\[
\langle\mathcal{L}u,u\rangle \geq \|\nabla u\|^2_{L^2(\R^2)}-C \|u\|^2_{L^2(\R^2)}.
\]
We conclude from  above two estimates that \eqref{1.76} holds true.

(2). If $a\in(0,\infty)$, it is obvious that $\mathcal{N}>\mathcal{L}\geq 0$ and thus $ker \mathcal{N}=\{0\}$.
On the other hand, if $a\in(-a^*,0)$, then the non-degenerate of $\mathcal{N}$ is a well-known result of \cite[Corollary 11.5]{BO}, which thus implies that \eqref{2.19} holds true. Moreover, the same argument of \cite[Theorem 5.3]{GLY} yields that $\mathcal{N}^{-1}$: $L^2(\R^2)\longmapsto L^2(\R^2)$ exists and is a continuous linear operator satisfying \eqref{2.20}.

(3). Similar to \cite[Theorem 5.3]{GLY}), since the operator $\mathcal{N}^{-1}$  commutes
with the angular momentum, one can get that if $a\neq 0$ and $\phi(|x|)\in L^2(\R^2)$ is radially symmetric, then $\psi(x)=\mathcal{N}^{-1}\phi(|x|)\in M(\R^2,\R)$ is also radially symmetric, and the proof is therefore complete. \qed

\vskip 0.05truein

Applying Lemma \ref{lem-2}, we next derive the following refined estimates.

\begin{lem}\label{prop-2.2}
Under the assumptions of Theorem \ref{thm1.2}, let $\{\alpha_n\}$ be a positive bounded sequence, and consider positive constants $C>0$ and $ K>0$, which are independent of $n$.
\begin{enumerate}
\item If $q_n$ satisfies
\begin{equation}\label{prop4.1:A}
|x^{\perp}\cdot\nabla  {q}_n |\leq C\alpha_n e^{-K|x|}\ \ \mbox{in}\ \,\R^2,
\end{equation}
then there exists a constant $C=C(a,K)>0$, depending only on $a$ and $K$, such that
\begin{equation}\label{prop4.1:B}
|r_n(x)|\leq C(a,K)\Omega_n\alpha_ne^{-K|x|},\ \, |\nabla r_n(x)|\leq C(a,K)\Omega_n\alpha_ne^{-\frac{K}{2}|x|} \,\ \hbox{in\,\ $\R^2$.}
\end{equation}

\item 	If $f_n(x)\in M(\R^2,\R)$ and $g_n(x)\in L^2(\R^2,\R)$ satisfy
\begin{equation}\label{4.1}
\mathcal{N}f_n(x)=g_n(x),\,\ |g_n(x)|\leq C\alpha_ne^{-K|x|}\,\ \hbox{in $\R^2$,}
\end{equation}
then there exists a constant $C(a,K)>0$, depending only on $a$ and $K$, such that
\begin{enumerate}
\item[(a)] If $a\neq 0$, then we have
\begin{equation}\label{4.2}
|f_n(x)|\leq C(a,K)\alpha_ne^{-K|x|},\,\ |\nabla f_n(x)|\leq C(a,K)\alpha_ne^{-\frac{K}{2}|x|}\,\ \hbox{in $\R^2$.}
\end{equation}
\item[(b)] If $a= 0$ and $f_n(x)$ satisfies $\inte f_n(x)u_0(x)dx=0$, then we have
\begin{equation}\label{4.2*}
|f_n(x)|\leq C(K)\alpha_ne^{-K|x|},\,\ |\nabla f_n(x)|\leq C(K)\alpha_ne^{-\frac{K}{2}|x|}\,\ \hbox{in $\R^2$.}
\end{equation}
\end{enumerate}
\end{enumerate}
\end{lem}

{\noindent \bf  Proof. }1. Recall from \eqref{rev-4:3} that $r_n$ satisfies the following equation
\begin{equation}\label{2.13}
(-\Delta+V(x)-\mu_n+a|u_n|^2)r_n=-\Omega_n(x^{\perp}\cdot \nabla q_n)\ \ \hbox{in}\,\ \R^2.
\end{equation}
Following \eqref{2.9}, \eqref{2.28M} and \eqref{1.76}, we then have
\[
\begin{aligned}
&-\Omega_n\inte \big[(x^{\perp}\cdot \nabla q_n)\,  r_n\big]dx\\
=&\inte \big[|\nabla r_n|^2+(V(x)-\mu_n+a|u_n|^2)r^2_n\big]dx\\
\geq&\inte (\mathcal{L}r_n\cdot r_n)dx+\inte \big[(\mu_0-\mu_n)+a(|u_n|^2-u_0^2)\big]r_n^2dx\\
\geq& \,\rho\|r_n\|^2_{H^1(\R^2)}-o(1)\inte r_n^2dx\geq \frac{\rho}{2}\|r_n\|^2_{H^1(\R^2)} .
\end{aligned}
\]
 We then get from \eqref{prop4.1:A} and above that
\[
\frac{\rho}{2}\|r_n\|^2_{H^1(\R^2)}\leq\Big|\Omega_n\inte \big[(x^{\perp}\cdot \nabla q_n)\, r_n\big]dx\Big|\leq C\Omega_n\alpha_n\|r_n\|_{L^2(\R^2)} ,
\]
which further implies that
\begin{equation}\label{2.14}
\|r_n\|_{H^1(\R^2)}\leq C\Omega_n\alpha_n,
\end{equation}
where $\alpha_n>0$ is as in (\ref{prop4.1:A}).
On the other hand, we derive from \eqref{2.13} that $ r_n^2$ satisfies the following equation
\[
\begin{aligned}
\Big[-\frac 12\Delta +\Big(V(x)-\mu _n+a|u_{n}|^2\Big) \Big]r_{n}^2+|\nabla r_n|^2=-\Omega_n (x^\bot\cdot\nabla q_n)\,r_n\ \ \hbox{in}\,\ \R^2,
\end{aligned}
\]
which then implies that
\begin{equation}\label{D20}
-\frac 12\Delta r_n^2-\mu _n r_n^2+a|u_n|^2 r_n^2\leq -\Omega_n (x^\bot\cdot\nabla q_n)r_n\ \ \hbox{in}\,\ \R^2.
\end{equation}
By De Giorgi-Nash-Moser theory \cite[Theorem 4.1]{HL}, we thus obtain from \eqref{D20} that for
any $y\in \R^2$,
\begin{equation}\label{D-21}
\sup_{x\in B_{\frac 12}(y)}|r_n(x)|^2\leq C\Big(\|r_n\|^2_{L^2(B_1(y))}+\|\Omega_n (x^\bot\cdot\nabla q_n)r_n\|_{L^2(B_1(y))}\Big).
\end{equation}
Combining \eqref{2.14} and \eqref{D-21} yields that
\begin{equation}\label{rev-28}
\|r_n\|_{L^\infty(\R^2)}\leq  C\Omega_n\alpha_n .
\end{equation}
Moreover, since $\mu _n\to \mu_0$ as $n\to \infty$, there exists a sufficiently large constant $R=R(a,K)>0$ such that
\[
V(x)-\mu_n+a|u_n|^2\geq K^2+1\,\ \hbox{in}\,\ \R^2/B_R(0).
\]
By the comparison principle, we then deduce from \eqref{prop4.1:A} and \eqref{2.13}  that
\[
|r_n(x)|\leq C\Omega_n\alpha_ne^{-K|x|}\,\ \hbox{in}\,\ \R^2/B_R(0),
\]
from which  and \eqref{rev-28} we derive that
\begin{equation}\label{2.17}
|r_n(x)|\leq C\Omega_n\alpha_n e^{-K|x|}\,\ \hbox{in}\,\ \R^2.
\end{equation}
Furthermore, applying gradient estimates (see (3.15) in \cite{GT}) to the equation \eqref{2.13}, we conclude from above that
\begin{equation}\label{D2}
|\nabla r_n(x)|\leq C\Omega_n\alpha_n e^{-\frac {K}{2}|x|}\ \ \hbox{in\,  $\R^2$},
\end{equation}
which therefore implies that (\ref{prop4.1:B}) holds true in view of (\ref{2.17}).

2. We first prove that \eqref{4.2} holds true, if $a\neq 0$. Since $f_n$ satisfies \eqref{4.1}, we derive from \eqref{2.20} that
\[
\begin{aligned}
\hat{\rho}\|f_n\|_{H^2(\R^2)}
\leq\|g_n\|_{L^2(\R^2)}\leq C(a,K)\alpha_n,
\end{aligned}
\]
which then implies that
\begin{equation}\label{4.3}
\|f_n\|_{L^\infty(\R^2)}\leq C(a,K)\alpha_n.
\end{equation}
By the comparison principle, we then derive from \eqref{4.1} that
\begin{equation}\label{4.4}
|f_n(x)|\leq C(a,K)\alpha_n e^{-K|x|}\ \ \hbox{in}\,\ \R^2.
\end{equation}
%Following \eqref{2.15} and \eqref{2.18}, there exists a large constant $R=R(\mu)>0$ such that
%\[
%\big|-\Delta w_n+2w_n\big|\leq C\e e^{-\frac {1}{16}|x|}\,\ \hbox{in}\,\ \R^2/B_R(0),
%\]
Applying again gradient estimates (see (3.15) in \cite{GT}), we further obtain from \eqref{4.4} that $|\nabla f_n|$ also satisfies the desired estimate (\ref{4.2}).

Next, we prove that \eqref{4.2*} holds true, if $a=0$ and  $f_n(x)$ satisfies $\inte f_n(x)u_0(x)dx=0$. Indeed, in this case, we have
\[
\mathcal{N}=\mathcal{L}=-\Delta+V(x)-\mu_0.
\]
Since $f_n$ satisfies \eqref{4.1} and $\inte f_n(x)u_0(x)dx=0$, we derive from  Lemma \ref{lem-2} (1) that
\begin{equation}\label{4.5}
\|f_n\|^2_{H^1(\R^2)}\leq C\langle \mathcal{L} f_n,f_n\rangle\leq C\|g_n\|_{L^2(\R^2)}\|f_n\|_{L^2(\R^2)}\leq C(K)\alp_n\|f_n\|_{L^2(\R^2)},
\end{equation}
which implies that
\begin{equation}\label{4.6}
\|f_n\|_{H^1(\R^2)}\leq C(K)\alp_n.
\end{equation}
Similar to the proof of \eqref{rev-28}, since $f_n$ satisfies \eqref{4.6}, we obtain from De Giorgi-Nash-Moser theory \cite[Theorem 4.1]{HL} that
\begin{equation}\label{4.7}
\|f_n\|_{L^\infty(\R^2)}\leq C(K)\alp_n.
\end{equation}
Following \eqref{4.7}, the comparison principle applied to \eqref{4.1} then yields that
\begin{equation}\label{4.8}
|f_n(x)|\leq C(K)\alpha_n e^{-K|x|}\ \ \hbox{in}\,\ \R^2,
\end{equation}
from which one can further obtain (\ref{4.2*}) on the similar exponential decay of $|\nabla f_n|$. This completes the proof of Lemma \ref{prop-2.2}.  \qed

\section{Proof of Theorem \ref{thm1.2}}

For any fixed $a\in[0,+\infty)$ or a.e. $a\in (-a^*,0)$, this section is devoted to the complete proof of Theorem \ref{thm1.2} on the nonexistence of vortices for minimizers of $e(\Omega,a)$ in the case where $\Omega>0$ is small enough, since the case $\Omega=0$ follows directly in view of Theorem B and those explanations below (\ref{2.27}). Under the assumptions of Theorem \ref{thm1.2}, consider the minimizer  $u_n$ of $e(\Omega,a)$, and let  $(q_n,r_n)$ be as in \eqref{4.10}. We also recall  from Theorem B that $u_0>0$ is the unique positive minimizer of $e(a)$, and $(u_0,\mu _0)$ satisfies the elliptic equation (\ref{3.4MM}).

Following the estimates of the previous section, we first establish the following corollary.

\begin{cor}\label{prop-1.1}
Under the assumptions of Theorem \ref{thm1.2}, define
\begin{equation}\label{4.9}
w_n(x):=q_n(x)-u_0(x),\quad \e:=\max\big\{\Omega_n,|\mu_n-\mu_0|\big\}>0,
\end{equation}
where $\mu_n$ and $\mu_0$ are as in (\ref{lem-2.1}) and (\ref{3.4MM}), respectively.
Then we have
\begin{enumerate}
\item There exists a constant $C>0$, independent of $n$,  such that
\begin{equation}\label{cor4.1:B}
|r_n(x)|\leq C\Omega_ne^{-|x|},\ \, |\nabla r_n(x)|\leq C\Omega_ne^{-\frac{1}{2}|x|}\ \ \mbox{in}\, \ \R^2.
\end{equation}
\item There exists a constant $C>0$, independent of $n$,  such that
\begin{enumerate}
\item[(a)] If $a\neq 0$, then $w_n$ satisfies
\begin{equation}\label{cor-1.2}
|w_n(x)|\leq C\e e^{-|x|},\ \, |\nabla w_n(x)|\leq C\e e^{-\frac{1}{2}|x|}\ \ \mbox{in}\, \ \R^2.
\end{equation}
\item[(b)] If $a=0$, then $w^{\perp}_n=w_n-\big({\inte w_nu_0dx}\big)u_0$ satisfies
\begin{equation}\label{cor-1.2*}
|w^{\perp}_n(x)|\leq C\e e^{-|x|},\ \, |\nabla w^{\perp}_n(x)|\leq C\e e^{-\frac{1}{2}|x|}\ \ \mbox{in}\, \ \R^2.
\end{equation}
\end{enumerate}
\end{enumerate}
\end{cor}

\noindent{\bf Proof.}  Following Lemma \ref{lem-2.1} (3), the estimates of \eqref{cor4.1:B} follow directly from  Lemma \ref{prop-2.2} by taking  $\alpha_n\equiv 1$ and $K=1$.

Next, we prove that \eqref{cor-1.2} holds if $a\neq 0$. Direct calculations give from (\ref{rev-4:3}) that $w_n$ satisfies the following equation
\begin{equation}\label{2.15}
\mathcal{N}w_n+a(q_n+2u_0)w_n^2=(\mu_n-\mu _0)q_n-ar_n^2q_n+\Omega_n(x^\perp\cdot \nabla r_n)\ \ \mbox{in}\, \ \R^2,
\end{equation}
where the linearized operator $\mathcal{N}$ is defined by (\ref{2.6}).
The left-hand side of \eqref{2.15} satisfies
\[
\begin{aligned}
\quad\|\mathcal{N}w_n+a(q_n+2u_0)w_n^2\|_{L^2(\R^2)}&\geq \|\mathcal{N}w_n\|_{L^2(\R^2)}-C(a)\|w_n\|^2_{L^4(\R^2)}\\
&\geq \hat{\rho}\|w_n\|_{H^2(\R^2)}-C(a)\|w_n\|^2_{H^2(\R^2)},
\end{aligned}
\]
where \eqref{2.20} is used in the second inequality. On the other hand, following \eqref{cor4.1:B}, the right-hand side of \eqref{2.15} satisfies
\[
\begin{aligned}
&\big|(\mu_n-\mu _0)q_n-ar_n^2q_n+\Omega_n(x^\perp\cdot \nabla r_n)\big|\leq C\e e^{-|x|}\ \ \hbox{in}\, \   \R^2.
\end{aligned}
\]
We thus derive from above that
\[
\hat{\rho}\|w_n\|_{H^2(\R^2)}-C(a)\|w_n\|^2_{H^2(\R^2)}\leq C\e,
\]
where Lemma \ref{lem-2} (2) is also used. This further implies that
\begin{equation}\label{2.16}
\|w_n\|_{L^\infty(\R^2)}\leq C\|w_n\|_{H^2(\R^2)}\leq C\e.
\end{equation}
Similar to the proofs of \eqref{rev-28} and \eqref{2.17}, one can further conclude from \eqref{2.15} and \eqref{2.16} that \eqref{cor-1.2} holds true.

Finally, we prove that (\ref{cor-1.2*}) holds for the case $a=0$. In this case, the operator $\mathcal{N}=\mathcal{L}=-\Delta+V(x)-\mu_0$, and $w_n^\perp=w_n-\big({\inte w_nu_0dx}\big)u_0$ satisfies $\inte w_n^\perp(x)u_0(x)dx=0$ and
\begin{equation}\label{4.11}
\mathcal{N}w_n^\perp=-\Delta w_n^\perp+V(x)w_n^\perp-\mu_0w_n^\perp=(\mu_n-\mu _0)q_n+\Omega_n(x^\perp\cdot \nabla r_n)\ \ \hbox{in}\,\ \R^2.
\end{equation}
Following Lemma \ref{lem-2.1} (3), we obtain that the right-hand side of \eqref{4.11} satisfies
\[
\big|(\mu_n-\mu _0)q_n+\Omega_n(x^\perp\cdot \nabla r_n)\big|\leq C\eps_ne^{-|x|}\,\ \hbox{in}\,\ \R^2.
\]
Thus, we derive from Lemma \ref{prop-2.2} (2) that \eqref{cor-1.2*} holds true, and we are done.  \qed

%\medskip
%{\noindent \bf Proof of Theorem \ref{thm1.2} in the case $a= 0$.} Note that in this case $\mathcal{N}=\mathcal{L}=-\Delta+V(x)-\mu_0$ and $w_n^{\perp}$ satisfies \eqref{4.11}.
%Following Corollary \ref{prop-1.1} (2b), we obtain that
%\[
%|x^{\perp}\cdot\nabla q_n|=|x^{\perp}\cdot\nabla w_n^\perp|\leq C\eps_ne^{-\frac{1}{4}|x|}.
%\]
%Therefore, we derive from Lemma \ref{prop-2.2} (1) that
%\begin{equation}\label{3.3}
%|r_n(x)|\leq C\Omega_n\e e^{-\frac 14|x|},\ \, |\nabla r_n(x)|\leq C\Omega_n\e e^{-\frac{1}{8}|x|}\ \ \mbox{in}\, \ \R^2.
%\end{equation}
%On the other hand, because $ker \mathcal{N}=ker \mathcal{L}=\{u_0\}$, we find
%\[
%\begin{aligned}
%&\quad\inte [(\mu_n-\mu _0)q_n+\Omega_n(x^\perp\cdot \nabla r_n)]u_0dx\\
%&=\inte (\mathcal{N}w_n^\perp u_0)dx
%=\inte (\mathcal{N}u_0w_n^\perp )dx=0.
%\end{aligned}
%\]
%Following \eqref{3.3}, we have $$\Big|\inte \Omega_n(x^\perp\cdot \nabla r_n)u_0dx\Big|\leq C\Omega_n^2\eps_n.$$ However, we derive from \eqref{2.3} that $$\Big|\inte (\mu_n-\mu _0)q_nu_0dx\Big|\geq C |\mu_n-\mu _0|,$$

 \vskip 0.05truein

\begin{lem}\label{lem4.5}
Under the assumptions of Theorem \ref{thm1.2}, for some $m\geq 2$ suppose $w_n=q_n-u_0$  can be decomposed as
\begin{equation}\label{D35}
w_n(x)=\psi_{m,n}(|x|)+T_{m,n}(x),
\end{equation}
where  $\psi_{m,n}(|x|)$ is radially symmetric, and $T_{m,n}$ satisfies
\begin{equation}\label{D33}
|T_{m,n}(x)|\leq A^m\e^{m}e^{-|x|},\ \, |\nabla T_{m,n}(x)|\leq A^m\e^{m}e^{-\frac{1}{2}|x|}\ \ \hbox{in\ \,$\R^2$}\, \ \mbox{as}\ \, n\to \infty,
\end{equation}
for some constant $A>0$, independent of $m$ and  $n$, where $\e>0$ is defined in \eqref{4.9}. Then the decomposition of (\ref{D35}) and (\ref{D33}) holds  for $m+1$.
\end{lem}

{\noindent \bf Proof.} Note that $q_n$ is the real part of $u_ne^{i\theta_n}$, $u_0>0$ and $\e >0$ are as in \eqref{4.9}. For some large constant $A>0$, we may assume that $\e >0$ is small enough so that $0<A\e^{\frac 12}<1$ for sufficiently large $n>0$. Moreover, we always use the same symbol $C$ to denote various positive constants which are independent of $A>0$, $m\ge 2$ and  $\e >0$.

Since $u_0(|x|)+\psi_{m,n}(|x|)$ is radially symmetric, we derive from \eqref{D35} and \eqref{D33} that for some $m\geq 2$,
\[
\big|x^{\perp}\cdot \nabla q_{n}\big|=\big|x^{\perp}\cdot \nabla T_{m,n}\big|\leq CA^{m}\e^{m}e^{-\frac{1}{4}|x|}\ \ \mbox{in\,\ $\R^2$},
\]
which and Lemma \ref{prop-2.2} (1) then yield that
\begin{equation}\label{G14}
|r_n(x)|\leq CA^{m}\e^{m+1}e^{-\frac{1}{4}|x|},\ \ |\nabla r_n(x)|\leq CA^{m}\e^{m+1}e^{-\frac{1}{8}|x|}\ \ \mbox{in\,\ $\R^2$}.
\end{equation}
The rest proof is divided into two steps.
\vskip 0.05truein

 {\em  Step 1.} We claim that there exists a decomposition $w_n=\psi_{m+1,n}+T_{m+1,n}$, where $\psi_{m+1,n}(x)=\psi_{m+1,n}(|x|)$ is radially symmetric and $T_{m+1,n}(x)$ satisfies
\begin{equation}\label{G17*}
|T_{m+1,n}|\leq CA^{m}\e^{m+1}e^{-\frac{1}{16}|x|},\,\ |\nabla T_{m+1,n}|\leq CA^{m}\e^{m+1}e^{-\frac{1}{32}|x|}\ \ \hbox{in\,\ $\R^2$.}
\end{equation}
To prove (\ref{G17*}), we shall consider separately the following two cases:
\vskip 0.05truein

%
%We shall prove the above claim by
%
%

{\em  Case 1:}
$a=0$. In this case, we define
\begin{equation}\label{G23M}
\psi_{m+1}\equiv \Big(\inte w_nu_0dx\Big)u_0,\ \ \hbox{and}\ \ T_{m+1}\equiv w_n^{\perp}=w_n-\Big(\inte w_nu_0dx\Big)u_0.
\end{equation}
Note from (\ref{4.11}) that for $a=0$, $\mathcal{N}=\mathcal{L}=-\Delta+V(x)-\mu_0$, $ker \mathcal{L}=\{u_0\}$, and $w_n^{\perp}$ satisfies
\begin{equation}\label{G23}
\mathcal{L}w_n^{\perp}=(\mu_n-\mu_0)q_n+\Omega_n(x^\perp\cdot\nabla r_n)\ \ \hbox{in}\ \ \R^2,\ \ \inte w_n^{\perp}u_0dx=0.
\end{equation}
Since
\[
0=\langle \mathcal{L}u_0, w_n^{\perp}\rangle=\langle u_0, \mathcal{L}w_n^{\perp}\rangle=\langle u_0, (\mu_n-\mu_0)q_n+\Omega_n(x^\perp\cdot\nabla r_n)\rangle,
\]
and $q_n\to u_0$ strongly in $L^\infty(\R^2,\R)$ as $n\to\infty$, we can deduce from \eqref{G14} that
\begin{equation}\label{G15}
\begin{split}
\frac{1}{2}|\mu_n-\mu_0|\langle u_0, u_0\rangle
&\leq |\langle u_0, (\mu_n-\mu_0)q_n\rangle|\\
&=|\langle u_0, \Omega_n(x^\perp\cdot\nabla r_n)\rangle|\leq CA^m\eps_n^{m+1}.
\end{split}
\end{equation}
%Then Lemma \ref{lem-2} (1) yields that
%\[
%\|w_n^{\perp}\|_{H^1(\R^2)}\leq CA^m\eps_n^{m+1},
%\]
%which further implies that
%\[
%\|w_n^{\perp}\|_{L^\infty(\R^2)}\leq CA^m\eps_n^{m+1}
%\]
%by De Giorgi-Nash-Moser theory \cite[Theorem 4.1]{HL}.
%Moreover, the same argument of proving Corollary \ref{prop-1.1} (2b) gives that
Therefore, we deduce from (\ref{2.25}), (\ref{G14}) and above that the right-hand side of \eqref{G23} satisfies
\[
|(\mu_n-\mu_0)q_n+\Omega_n(x^\perp\cdot\nabla r_n)|\leq CA^m\eps_n^{m+1}e^{-\frac{1}{16}|x|}\,\ \hbox{in\,\ $\R^2$}.
\]
By Lemma \ref{prop-2.2} (2b), we thus obtain from (\ref{G23}) that
\[
|w_n^{\perp}(x)|\leq CA^m\eps_n^{m+1}e^{-\frac{1}{16}|x|},\,\ |\nabla w_n^{\perp}(x)|\leq CA^m\eps_n^{m+1}e^{-\frac{1}{32}|x|}\ \ \hbox{in}\, \ \R^2,
\]
which shows that the claim \eqref{G17*} holds true in the case $a=0$.
\vskip 0.05truein

 {\em  Case 2:} $a\neq 0$. In this case, for above $m\geq 2$  we obtain from (\ref{cor-1.2}), \eqref{D35} and \eqref{D33} that
\begin{equation}\label{D40}
\begin{split}
|\psi_{m,n}(x)|\leq A^m\e^{m}e^{-|x|}+C\e e^{-|x|}
\leq C\e e^{-|x|} \quad\hbox{in}\,\ \R^2.
\end{split}
\end{equation}
Similarly, we have
\begin{equation}\label{D41}
|\nabla\psi_{m,n}(x)|\leq C\e e^{-\frac{1}{2}|x|} \quad\hbox{in}\,\ \R^2.
\end{equation}
Since $q_n=u_0+\psi_{m,n}+T_{m,n}$,  we get from \eqref{2.15} that
\begin{equation}\label{G1}
\begin{split}
\mathcal{N}w_n:&=(\mu_n-\mu _0)(u_0+w_n)-a(q_n+2u_0)w_n^2-ar_n^2q_n+\Omega_n(x^\perp\cdot \nabla r_n)\\
&=\Omega_n\big(x^{\perp}\cdot\nabla r_n\big)-ar_n^2q_n+B_{1,n}(x)+B_{2,n}(x)\quad\hbox{in}\,\ \R^2,
\end{split}
\end{equation}
where the radially symmetric function $B_{1,n}(x)$ satisfies
\begin{equation}\label{G4}
\begin{split}
B_{1,n}(x)&=(\mu_n-\mu _0)(u_0+\psi_{m,n})-a\big(\psi_{m,n}^3+3u_0\psi_{m,n}^2\big),
\end{split}
\end{equation}
and the non-radially symmetric function $B_{2,n}(x)$ satisfies
\begin{equation}\label{G5}
\begin{split}
B_{2,n}(x)=&(\mu_n-\mu _0)T_{m,n}-a\big(3\psi_{m,n}^2+3\psi_{m,n}T_{m,n}\\
&+T^2_{m,n}+6u_0\psi_{m,n}+3u_0T_{m,n}\big)T_{m,n}.
\end{split}
\end{equation}
Following Lemma \ref{lem-2} (3) and (\ref{G4}), there exists a unique function $\psi_{m+1,n}\in C^2(\R^2)\cap L^\infty(\R^2)$ of
\begin{equation}\label{G6}
\begin{split}
\mathcal{N}\psi_{m+1,n}=B_{1,n}(x)\quad\hbox{in\  $\R^2$,}
\end{split}
\end{equation}
and moreover, $\psi_{m+1,n}$ is also radially symmetric.

We now define
\begin{equation}\label{G7}
T_{m+1,n}(x):=w_n(x)-\psi_{m+1,n}(x),
\end{equation}
where $\psi_{m+1,n}(x)$ is given by (\ref{G6}), so that $T_{m+1,n}(x)$ satisfies
\begin{equation}\label{G16}
\mathcal{N}T_{m+1,n}=\Omega_n\big(x^{\perp}\cdot\nabla r_n\big)-ar_n^2q_n+B_{2,n}(x)
\ \ \mbox{in}\,\  \R^2,
\end{equation}
due to (\ref{G1}).
%Therefore, under the decomposition of (\ref{G7}), the rest is to prove that (\ref{D33}) holds  for $m+1$.
%In the following, we obtain rough estimates of $T_{m+1,n}$ and $\nabla T_{m+1,n}$.
Following \eqref{D33}, \eqref{G14} and \eqref{D40}, we deduce that there exists a constant $C>0$ such that
\begin{equation}\label{5.1}
|\Omega_n(x^\perp\cdot\nabla r_n)-ar_n^2q_n|\leq CA^m\e^{m+1}e^{-\frac{1}{16}|x|}\ \ \hbox{in}\, \ \R^2,
\end{equation}
and
\begin{equation}\label{5:G14}
|B_{2,n}(x)|\leq CA^m\e^{m+1}e^{-|x|},\,\ \ |\nabla B_{2,n}(x)|\leq CA^m\e^{m+1}e^{-\frac{1}{2}|x|}\ \ \hbox{in}\, \ \R^2.
\end{equation}
In view of (\ref{5.1}) and (\ref{5:G14}), we conclude from Lemma \ref{prop-2.2} (2$a$) that the claim \eqref{G17*} also holds true for the case $a\neq 0$.

%Moreover, we infer from the exponential decay of $\nabla I_a, I_a,B_{2,a}$ that
%\[
%|\mathcal{N}_aT_{m+1,a}|\leq CA^m\e^{2(m+1)}e^{-\frac{1}{16}|x|}\,\,\hbox{in $\R^2/B_R(0)$.}
%\]
%We then deduce from the estimate \eqref{G17*} and the comparison principle that
%\[
%|T_{m+1,a}|\leq CA^m\e^{2(m+1)}e^{-\frac{1}{16}|x|}\,\,\hbox{in $\R^2/B_R(0)$.}
%\]
%The estimate \eqref{G17*} and the above decay of $T_{m+1,a}$ further implies that
%\[
%\|\nabla T_{m+1,a}\|_{L^{\infty}(\R^2)}\leq CA^{m}\e^{2(m+1)},
%\]
%and therefore we conclude that
%\begin{equation}\label{G17}
%\| T_{m+1,a}\|_{C^{1}(\R^2)}\leq CA^{m}\e^{2(m+1)}.
%\end{equation}

\vskip 0.05truein
\noindent{\em  Step 2.} Under the decompositions (\ref{G23M}) or \eqref{G7}, we shall prove that (\ref{D33}) holds for $m+1$, which then completes the proof of the lemma. Towards this aim, actually it suffices to show that $T_{m+1,n}$ decays faster than \eqref{G17*} as $|x|\to\infty$.

Recall from (\ref{rev-4:3}), (\ref{G23}), (\ref{G7}) and (\ref{G16}) that $(r_n,T_{m+1,n})$ satisfies the following system
\begin{equation}\label{G8}
\left\{
\begin{aligned}
\mathcal{N}T_{m+1,n}&=\Omega_n\big(x^{\perp}\cdot\nabla r_n\big)-ar_n^2q_n+\hat B_{2,n}(x)
\ \ &\mbox{in} \ \ \R^2,\\
\mathcal{L}_n r_n&=-\Omega_n (x^\bot\cdot\nabla T_{m+1,n})
\ \ &\mbox{in} \,\ \R^2,
\end{aligned}
\right.
\end{equation}
where $\hat{B}_{2,n}(x)=(\mu_n-\mu_0)q_n$ if $a=0$, and $\hat{B}_{2,n}(x)={B}_{2,n}(x)$ if $a\neq 0$.
Multiplying the first equation of \eqref{G8} by $T_{m+1,n}$ and the second one of \eqref{G8} by $r_n$, respectively, we then have
\[
\left\{
\begin{aligned}
&
-\frac{1}{2}\Delta |T_{m+1,n}|^2+|\nabla T_{m+1,n}|^2+\Big[V(x)-\mu_0+3au_0^2\Big]|T_{m+1,n}|^2\\
&\quad\quad=\Big[\Omega_n\big(x^{\perp}\cdot\nabla r_n\big)-ar_n^2q_n+\hat B_{2,n}\Big]T_{m+1,n}\quad\mbox{in\,\ $\R^2$},\\
&-\frac{1}{2}\Delta  r_n^2+|\nabla r_n|^2+\Big[V(x)-\mu_n+a|u_n|^2\Big] r_n^2=-\Omega_n (x^\bot\cdot\nabla T_{m+1,n})r_n\quad\mbox{in\,\ $\R^2$}.
\end{aligned}
\right.
\]
By the diamagnetic inequality,  we have
\begin{equation}\label{5:A:M}
|\nabla T_{m+1,n}|^2+\frac{\Omega^2_n|x|^2}{4}| r_n|^2+\Omega_n (x^\bot\cdot\nabla T_{m+1,n})r_n\geq 0\quad\mbox{in\,\ $\R^2$},
\end{equation}
and
\[
|\nabla r_n|^2+\frac{\Omega^2_n|x|^2}{4}|T_{m+1,n}|^2-\Omega_n (x^\bot\cdot\nabla r_n)T_{m+1,n}\geq 0\quad\mbox{in\,\ $\R^2$}.
\]
By the exponential decay (\ref{2.25}) of $|u_n|$, we also get from  \eqref{G14}, (\ref{G15}) and \eqref{5:G14} that as $n\to \infty$,
\[
\begin{aligned}
\big(-ar_n^2q_n+\hat B_{2,n}\big)T_{m+1,n}&\leq \frac{1}{2}|T_{m+1,n}|^2+\frac{1}{2}\big(-ar_n^2q_n+\hat B_{2,n}\big)^2\\
&\leq \frac{1}{2}|T_{m+1,n}|^2+C\e^{2(m+1)}A^{2m}e^{-2|x|}\quad\mbox{in\,\ $\R^2/B_R$},
\end{aligned}
\]
where and below $R>0$ is sufficiently large.
We then obtain from above that as $n\to \infty$,
\[
\frac{1}{2}V(x)-\mu_0+3a|u_0|^2+a|u_n|^2\geq 3\quad\hbox{in}\,\ \R^2/B_R ,
\]
which further implies that
\begin{equation}\label{G18}
\begin{split}
\quad&-\frac{1}{2}\Delta\big(|T_{m+1,n}|^2+ r_n^2\big)+\frac{5}{2}\big(|T_{m+1,n}|^2+ r_n^2\big)\\
\leq &\,CA^{2m}\e^{2(m+1)}e^{-2|x|}\quad\mbox{in\ \ $\R^2/B_R$}.
\end{split}
\end{equation}
By the comparison principle, we thus derive from \eqref{G17*} and \eqref{G18} that as $n\to\infty$,
\begin{equation}\label{G19}
|T_{m+1,n}(x)|,\ |r_n(x)|\leq CA^m\e^{m+1}e^{-|x|}\quad\hbox{in\ \ $\R^2$.}
\end{equation}

We next address the desired exponential decay of $|\nabla T_{m+1,n}|$ and $|\nabla r_{n}|$.
Note from \eqref{G8} that $\partial_iT_{m+1,n}:=\frac{\partial T_{m+1,n}}{\partial x_i}$ satisfies
\[
\begin{aligned}
&\quad\mathcal{N}\partial_iT_{m+1,n}+\Big[\partial_iV(x)+6au_0\partial_iu_0\Big]T_{m+1,n}\\
&=\partial_i\Big[\Omega_n\big(x^{\perp}\cdot\nabla r_n\big)-ar_n^2q_n+\hat B_{2,n}(x)\Big]\quad\mbox{in\,\ $\R^2$}, \ i=1,\,2.
\end{aligned}
\]
Multiplying the above equation by $\partial_iT_{m+1,n}$ and summing it for $i=1,2$, we obtain that
\begin{equation}\label{G20}
\begin{split}
&\quad-\frac{1}{2}\Delta|\nabla T_{m+1,n}|^2+\Big[V(x)-\mu_0+3au_0^2\Big]|\nabla T_{m+1,n}|^2\\
&\quad +\sum_{i=1}^2|\nabla \partial_iT_{m+1,n}|^2+\sum_{i=1}^2\Big[\partial_iV(x)+6au_0\partial_iu_0\Big]T_{m+1,n}\partial_iT_{m+1,n}\\
&=\sum_{i=1}^2\partial_i\Big[\Omega_n\big(x^{\perp}\cdot\nabla r_n\big)-ar_n^2q_n+\hat B_{2,n}(x)\Big]\partial_iT_{m+1,n}\quad\mbox{in\,\  $\R^2$}.
\end{split}
\end{equation}
By the exponential decay   of $u_0$ and $|u_n|$, we calculate from (\ref{G17*}) and (\ref{G19}) that as $n\to \infty$,
\begin{equation}\label{G20M}
\begin{aligned}
&\quad\sum_{i=1}^2\Big[\partial_iV(x)+6au_0\partial_iu_0\Big]T_{m+1,n}\partial_iT_{m+1,n}\\
&\leq \sum_{i=1}^2\Big[|\partial_iV(x)|^2+C\Big]|T_{m+1,n}|^2 +\sum_{i=1}^2\frac{1}{4}|\partial_iT_{m+1,n}|^2\\
&\leq CA^{2m}\e^{2(m+1)}e^{-|x|}+\frac{1}{4}|\nabla T_{m+1,n}|^2 \quad \hbox{in\,\ $\R^2/B_R$},
\end{aligned}
\end{equation}
where and below $R>0$ is as before sufficiently large, and the assumption $(V)$ is used in the last inequality.

Similar to (\ref{G20M}), we have
\[
\begin{aligned}
&\quad\sum_{i=1}^2\partial_i\Big[-ar_n^2q_n+\hat B_{2,n}(x)\Big]\partial_iT_{m+1,n}\\
&\leq \frac{1}{4}\sum_{i=1}^2|\partial_iT_{m+1,n}|^2+\sum_{i=1}^2\Big[\partial_i\big(-ar_n^2q_n+\hat B_{2,n}\big)\Big]^2\\
&\leq \frac{1}{4}|\nabla T_{m+1,n}|^2+CA^{2m}\e^{2(m+1)}e^{-|x|}\quad \mbox{in\,\ $\R^2/B_R$,}
\end{aligned}
\]
which further implies that as $n\to\infty$,
\begin{equation}\label{G21}
\begin{split}
&\quad-\frac{1}{2}\Delta|\nabla T_{m+1,n}|^2+\Big[V(x)-\mu_0+3au_0^2-\frac 12\Big]|\nabla T_{m+1,n}|^2+\sum_{i=1}^2|\nabla \partial_iT_{m+1,n}|^2\\
&\leq\sum_{i=1}^2\partial_i\Big[\Omega_n\big(x^{\perp}\cdot\nabla r_n\big)\Big]\partial_iT_{m+1,n}+CA^{2m}\e^{2(m+1)}e^{-|x|}\quad\mbox{in\,\ $\R^2/B_R$},
\end{split}
\end{equation}
in view of (\ref{G20}) and (\ref{G20M}).
Similar to \eqref{G21}, we also get that as $n\to \infty$,
\begin{equation}\label{G22}
\begin{split}
&\quad-\frac{1}{2}\Delta |\nabla r_n|^2+\Big[V(x)-\mu_n+a|u_n|^2-\frac 12\Big]|\nabla r_n|^2+\sum_{i=1}^2|\nabla \partial_ir_n|^2\\
&\leq-\sum_{i=1}^2\partial_i\Big[\Omega_n\big(x^{\perp}\cdot\nabla T_{m+1,n}\big)\Big]\partial_ir_n+CA^{2m}\e^{2(m+1)}e^{-|x|}\quad\mbox{in\,\ $\R^2/B_R$}.
\end{split}
\end{equation}
Using the diamagnetic inequality as in (\ref{5:A:M}), one can derive from \eqref{G21} and \eqref{G22} that as $n\to \infty$,
\[
-\frac{1}{2}\Delta\Big(|\nabla r_n|^2+|\nabla T_{m+1,n}|^2\Big)+\Big(|\nabla r_n|^2+|\nabla T_{m+1,n}|^2\Big)\leq CA^{2m}\e^{2(m+1)}e^{-|x|}\quad \mbox{in\,\ $\R^2/B_R$.}
\]
By the comparison principle, we then derive from above that as $n\to \infty$,
\begin{equation}\label{G29}
|\nabla r_n|,\ \ |\nabla T_{m+1,n}|\leq CA^{m}\e^{(m+1)}e^{-\frac 12|x|}\quad \mbox{in\ \ $\R^2/B_R$.}
\end{equation}
Therefore, we conclude from \eqref{G14}, \eqref{G17*} and \eqref{G29} that there exists a sufficiently large $C>0$ such that  as $n\to \infty$,
\begin{equation}\label{5:G29}
|\nabla r_n(x)|,\ \ |\nabla T_{m+1,n}(x)|\leq CA^{m}\e^{(m+1)}e^{-\frac{1}{2}|x|}\quad \mbox{in\,\ $\R^2$.}
\end{equation}

Because the positive constant $C$ in (\ref{5:G29}) is independent of the values $A>0$, $m\ge 2$ and  $n$, one can choose a sufficiently large constant $A$ such that $A>C$. We thus conclude that (\ref{D33}) holds for $m+1$ in view of \eqref{G19} and (\ref{5:G29}), and the proof of Lemma \ref{lem4.5} is therefore complete. \qed

\vskip 0.05truein

{\noindent \bf Proof of Theorem \ref{thm1.2}.}
Set $T_{1,n}(x)=w^\perp_n(x)=w_n-\big(\inte w_nu_0dx\big)u_0$ and $\psi_{1,n}=\big(\inte w_nu_0dx\big)u_0$ if $a=0$,
while $T_{1,n}(x)=w_n(x)$ and $\psi_{1,n}\equiv 0$ if $a\neq 0$. We then obtain from Corollary \ref{prop-1.1} (2) that
\begin{equation}\label{5A:G29}
|T_{1,n}(x)|\leq C_1\e e^{-|x|},\ \ |\nabla T_{1,n}(x)|\leq C_1\e e^{-\frac 12 |x|}\ \ \mbox{in}\ \, \R^2,
\end{equation}
where $C_1>0$ is independent of  $n>0$.
Set $T_{2,n}(x)=T_{1,n}(x)$ if $a=0$, and $T_{2,n}:=q_n-u_0-\psi_{2,n}$ if $a\neq 0$, where
$\psi_{2,n}\in C^2(\R^2)\cap L^\infty(\R^2)$ is the unique solution of
\[
\mathcal{N}
\psi_{2,n}= (\mu_n-\mu _0)u_0\ \ \mbox{in}\ \, \R^2.
\]
It then follows from Lemma \ref{lem-2} (3) that $\psi_{2,n}(|x|)$ is radially symmetric. Moreover, based on (\ref{5A:G29}), the same argument of proving Lemma \ref{lem4.5} gives that there exists a constant $C_2>0$, independent of $n>0$, such  that
\[
|T_{2,n}(x)|\leq C_2\e^2e^{-|x|},\ \ |\nabla T_{2,n}(x)|\leq C_2\e^2e^{-\frac 12 |x|}\ \ \mbox{in}\ \, \R^2\quad \mbox{as}\ \ n\to \infty.
\]
Take $A>0$ large enough that $A^2>C_2$, from which we have
\begin{equation}\label{5B:G29}
|T_{2,n}(x)|\leq A^2\e^2e^{-|x|},\,\  |\nabla T_{2,n}(x)|\leq A^2\e^2e^{-\frac 12 |x|}\ \ \mbox{in}\ \, \R^2\quad \mbox{as}\ \ n\to \infty.
\end{equation}

By Lemma \ref{lem4.5}, we deduce from (\ref{5B:G29}) that for any $m\geq 2$,
\[
| T_{m,n}(x)|\leq A^{m}\e^{m}e^{-|x|},\ \ |\nabla T_{m,n}(x)|\leq A^{m}\e^{m}e^{-\frac 12|x|}\ \ \mbox{in}\ \, \R^2\quad \mbox{as}\ \ n\to \infty.
\]
Recall from (\ref{D35}) that $q_n(x)=[u_0(x)+\psi_{m,n}(|x|)]+T_{m,n}(x)$, where $u_0(x)+\psi_{m,n}(|x|)$ is radially symmetric.
Applying Lemma \ref{prop-2.2}, we thus derive from above that for any $m\geq 2$,
\begin{equation}\label{5C:G29}
|r_n(x)|\leq CA^m\e^{m+1}e^{-\frac 14|x|},\,\ |\nabla r_n(x)|\leq CA^m\e^{m+1}e^{-\frac 18|x|}\ \ \mbox{in}\ \, \R^2\quad \mbox{as}\,\ n\to \infty.
\end{equation}
Therefore, we conclude from (\ref{5C:G29}) that for any $m\geq 2$,
\[
\|r_n\|_{C^1(\R^2)}\leq C\e^{\frac{m+2}{2}}\ \ \mbox{in}\ \,\R^2 \,\  \mbox{for sufficiently large}\,\ n,
\]
which further implies that $r_n\equiv 0$ as $n\to\infty$.

Since $r_n(x)\equiv 0$ as $n\to \infty$, $q_n$ must be a positive minimizer of \eqref{2.2}, and hence $q_n$ satisfies
\begin{equation}\label{5D:G29}
-\Delta q_n+\big(V(x)-\mu_0+aq_n^2\big)q_n=0\ \ \mbox{in}\,\ \R^2,
\end{equation}
which implies the absence of vortices for $q_n$ as $n \to \infty$.	
Following Theorem B, we also obtain the uniqueness of $q_n\equiv u_0$ for  any fixed $a\in[0,+\infty)$ or a.e. $a\in (-a^*,0)$. This therefore completes the proof of Theorem \ref{thm1.2}. \qed


\vskip 0.05truein




\begin{thebibliography}{GNN}
	
\bibitem{Abo} J. R. Abo-Shaeer, C. Raman, J. M. Vogels and W. Ketterle, {\em Observation of vortex lattices in Bose-Einstein condensate}, Science {\bf 292} (2001), 476.
	
	
	
	%\bibitem{Adamas} R. A. Adamas, Sobolev Spaces, Pure and Applied Mathematics, 65, Academic Press., New York, 1975.
	
	
	
	
	
\bibitem{A} A. Aftalion,  Vortices in Bose-Einstein condensates, Progress in Nonlinear Differential Equations and their Applications  67, Birkh$\ddot{a}$user Boston, Inc., Boston, MA,  2006.

\bibitem{AB}  A. Aftalion, X. Blanc and F. Nier, {\em Lowest Landau level functional and Bargmann spaces for Bose-Einstein condensates}, J. Funct. Anal. {\bf 241} (2006), 661--702.
	

	
	
	
%\bibitem{AA} A. Aftalion, S. Alama and L. Bronsard, {\em Giant vortex and the breakdown of strong pinning in a rotating Bose-Einstein condensate}, Arch. Rational Mech. Anal. {\bf 178} (2005), 247--286.
	%
	%
	%
	%
	%
\bibitem{AJ} A. Aftalion, R. L. Jerrard and J. Royo-Letelier, {\em Non-existence of vortices in the small density region of a condensate}, J. Funct. Anal. {\bf 260} (2011), 2387--2406.
	
\bibitem{AN} A. Aftalion, B. Noris, and C. Sourdis, {\em Thomas-Fermi approximation for coexisting two component Bose-Einstein condensates and nonexistence of vortices for small rotation}, Comm. Math. Phys. {\bf 336} (2015), 509--579.

	
	
\bibitem{Al} G. Allaire, {\em Homogenization and two-scale convergence}, SIAM J. Math. Anal. {\bf 23} (1992), 1482--1518.	
	
	
	
%\bibitem{Anderson} M. H. Anderson, J. R. Ensher, M. R. Matthews, C. E. Wieman and E. A. Cornell, {\em Observation of Bose-Einstein condensation in a dilute atomic vapor}, Science {\bf 269} (1995), 198--201.
%	
%	
%	
%\bibitem{ANS} J. Arbunich, I. Nenciu and C. Sparber, {\em Stability and instability properties of rotating Bose-Einstein condensates}, Lett. Math. Phys. {\bf 109} (2019), 1415--1432.
%	
	
	
	
	
	
	
	
	
%\bibitem{AS} G. Arioli and A. Szulkin, {\em A semilinear Schr$\ddot{o}$dinger equation in the presence of a magnetic field}, Arch. Ration. Mech. Anal. {\bf 170} (2003), 277--295.
	
	
%	
%\bibitem{BC} W. Bao and Y. Cai, {\em Ground states of two-component Bose-Einstein condensates with an internal atomic Josephson junction}, East Asia J. Appl. Math. {\bf 1} (2011), 49--81.
	
%\bibitem{Peng} T. Bartsch, E. N. Dancer and S. Peng. {\em On multi-bump semi-classical bound states of nonlinear Schr\"{o}dinger equations with electromagnetic fields}, Adv. Differential Equations {\bf 11} (2006), 781--812.
	
	
	
%\bibitem{B} I. Bloch, J. Dalibard and W. Zwerger, {\em Many-body physics with ultracold gases}, Reviews of Modern Phys. {\bf 80} (2008), 885--964.
%	
%	
%	
%	
%	
%	
%	
%\bibitem{Hulet2} C. C. Bradley, C. A. Sackett and R. G. Hulet, {\it Bose-Einstein condensation of lithium: observation of limited condensate number}, Phys. Rev. Lett. {\bf 78} (1997), 985.
%	
%	
%	
%\bibitem{Hulet1} C. C. Bradley, C. A. Sackett, J. J. Tollett and R. G. Hulet, {\it Evidence of Bose-Einstein condensation in an atomic gas with attractive interactions}, Phys. Rev. Lett. {\bf 75} (1995), 1687. {\it Erratum} Phys. Rev. Lett. {\bf 79} (1997), 1170.
%	
	
	
	
	
\bibitem{BO} J. Byeon and Y. Oshita, {\em Uniqueness of standing waves for nonlinear Schr$\ddot{o}$dinger equations}, Proc. Roy. Soc. Edinburgh Sect. A {\bf  138} (2008), 975--987.
	%
	%
	%
	%\bibitem{BWang} J. Byeon and Z. Q. Wang, {\em Standing waves with a critical frequency for nonlinear Schr$\ddot{o}$dinger equations}, Arch. Rational Mech. Anal. {\bf 165} (2002), 295--316.
	
%\bibitem{CLL} D. Cao, S. Li and P. Luo, {\em Uniqueness of positive bound states with multibump
%		for nonlinear Schrödinger equations}, Calc. Var. Partial Differential Equations
%	{\bf 54} (2015), 4037--4063.
	
%\bibitem{Cao} D. Cao and  Z. Tang, {\em Existence and uniqueness of multi-bump bound states of nonlinear  Schrodinger  equations with electromagnetic fields},  J. Differential Equations  {\bf 222}  (2006), 381--424.
	
	
%\bibitem{CC} L. D. Carr and C. W. Clark, {\em Vortices in attractive Bose-Einstein condensates in two dimensions}, Phys.
%	Rev. Lett. {\bf 97} (2006), 010403.
%	
%	
%\bibitem{C} T. Cazenave, Semilinear Schr$\ddot{o}$dinger equations, Courant Lecture Notes in Mathematics  Vol. 10, Courant Institute of Mathematical Science/AMS, New York, 2003.
	
\bibitem{CK} K. Chandrasekharan, Elliptic functions,  Grundlehren Math. Wiss., Springer, Berlin, 1985.
	
	
%\bibitem{CO} N. R. Cooper, {\em Rapidly rotating atomic gases,} Adv. Phys. {\bf 57} (2008),  539--616.
	
	
	
	
	%\bibitem{CD}  M. Correggi and D. Dimonte, {\em On the third critical speed for rotating Bose-Einstein condensates}, J. Math. Phys. {\bf 57} (2016), 071901.
	%
	%
	%
	%\bibitem{CP} M. Correggi, F. Pinsker, N. Rougerie and J. Yngvason, {\em Critical rotational speeds for superfluids in homogeneous traps}, J. Math. Phys. {\bf 53} (2012), 095203.
	
	
	
%\bibitem{CR} M. Correggi and N. Rougerie, {\em Boundary behavior of the Ginzburg-Landau order parameter in the surface superconductivity regime}, Arch. Rational Mech. Anal. {\bf 219} (2016), 553--606.
	
	
	
\bibitem{CRY} M. Correggi, N. Rougerie and J. Yngvason, {\em The transition to a giant vortex phase in a fast rotating Bose-Einstein condensate}, Comm. Math. Phys. {\bf 303} (2011), 451--508.
	
	
\bibitem{D} F. Dalfovo, S. Giorgini, L. P. Pitaevskii and S. Stringari, {\em Theory of Bose-Einstein condensation in trapped gases},  Rev. Math. Phys. {\bf 71} (1999), 463--512.
	
%\bibitem{Deng} Y. Deng, C. Lin and S. Yan, {\em On the prescribed scalar curvature problem in $\R^N$,
%		local uniqueness and periodicity}, J. Math. Pures Appl. {\bf 104} (2015), 1013--1044.
	
\bibitem{EL} M. J. Esteban and P. L. Lions, {\em Stationary solutions of nonlinear Schr$\ddot{o}$dinger equations with an external magnetic field}, Partial differential equations and the calculus of variations, Vol. I,  401--449, Progr. Nonlinear Differential Equations Appl. 1, Birkhuser Boston, Boston, MA, 1989.
	
	
\bibitem{F} A. L. Fetter, {\em Rotating trapped Bose-Einstein condensates},  Rev. Math. Phys. {\bf 81} (2009), 647--691.
	
	
	
\bibitem{GNN} B. Gidas, W. Ni and L. Nirenberg, {\em Symmetry of positive solutions of nonlinear elliptic equations in $\R^n$}, Mathematical analysis and applications  Part A, Adv. in Math. Suppl. Stud. Vol. {\bf 7}, Academic Press, New York  (1981), 369--402.
	
	
\bibitem{GT} D. Gilbarg and N. S. Trudinger, Elliptic Partial Differential Equations of Second Order, Springer, 1997.
	
	
%\bibitem{G60} E. P. Gross, {\em Structure of a quantized vortex in boson systems}, Nuovo Cimento {\bf 20} (1961), 454--466.
%	
%	
%\bibitem{G63} E. P. Gross, {\em Hydrodynamics of a superfluid condensate}, J. Math. Phys. {\bf 4} (1963), 195--207.
	
%\bibitem{Grossi} M. Grossi, {\em On the number of single-peak solutions of the nonlinear Schr\"odinger equations}, Ann. Inst. H. Poincar\'e Anal. Non Lin\'eaire  {\bf 19} (2002), 261--280.
	
	%\bibitem{GZ} H. L. Guo and H. S. Zhou, {\em A constrained variational problem arising in attractive Bose–Einstein condensate with ellipse-shaped potential}, Appl. Math. Lett. {\bf 87} (2019), 35--41.
	
%\bibitem{GL} Y. Guo, S. Li, J. Wei and X. Zeng, {\em Ground states of two-component attractive Bose-Einstein condensates II: semi-trivial limit behavior}, Trans. Amer. Math. Soc. {\bf 371} (2019), 6903--6948.
	
	
\bibitem{GLW} Y. Guo, C. Lin and J. Wei, {\em Local uniqueness and refined spike profiles of ground states for two-dimensional attractive Bose-Einstein condensates}, SIAM J. Math. Anal. {\bf 49} (2017), 3671--3715.


\bibitem{GLP} Y. Guo, Y. Luo and S. Peng, {\em Local uniqueness of ground states for rotating Bose-Einstein Condensates with attractive interactions}, submitted (2020), 28 pages, arXiv:2009.08013
	
\bibitem{GLY} Y. Guo, Y. Luo and W. Yang, {\em The nonexistence of vortices for  rotating Bose-Einstein condensates with attractive interactions},  Arch. Rational Mech. Anal. {\bf 238} (2020), 1231--1281.
	
	
	
	
	
	
\bibitem{GS}  Y. Guo and R. Seiringer, {\em On the mass concentration for Bose-Einstein condensates with attractive interactions}, Lett. Math. Phys. {\bf 104} (2014), 141--156.
	
	
	
\bibitem{GWZZ} Y. Guo, Z. Wang, X. Zeng and H. Zhou,  {\em  Properties of ground states of attractive Gross-Pitaevskii equations with multi-well potentials}, Nonlinearity {\bf 31}  (2018),  957--979.
	
	
	
%\bibitem{GZZ} Y. Guo, X. Zeng and H. Zhou, {\em Energy estimates and symmetry breaking  in attractive Bose-Einstein condensates with ring-shaped potentials}, Ann. Inst. H. Poincar\'e Anal. Non Lin\'eaire {\bf 33} (2016), 809--828.
	
	
	
\bibitem{HL} Q. Han and F. Lin,  Elliptic Partial Differential Equations, Courant Lecture Note in Math. 1, Courant Institute of Mathematical Science/AMS, New York, 2011.
	
	
	
%\bibitem{HM}  C. Huepe, S. Metens, G. Dewel, P. Borckmans and  M.E. Brachet, {\em Decay rates in attractive Bose-Einstein condensates}, Phys. Rev. Lett. {\bf 82}  (1999), 1616--1619.
	
	
	
\bibitem{IM-1} R. Ignat and V. Millot, {\em The critical velocity for vortex existence in a two-dimensional rotating Bose-Einstein condensate}, J. Funct. Anal. {\bf 233} (2006), 260--306.
	
	
	
\bibitem{IM-2}   R. Ignat and V. Millot, {\em Energy expansion and vortex location for a two-dimensional rotating Bose-Einstein condensate}, Rev. Math. Phys. {\bf 18} (2006), 119--162.
	
	
	
	
	
%\bibitem{KM} Y. Kagan, A. E. Muryshev and G.V. Shlyapnikov, {\em Collapse and Bose-Einstein condensation in a trapped Bose gas with nagative scattering length}, Phys. Rev. Lett. {\bf 81} (1998), 933--937.
%	
%	
%	
%
%
%	
%	
%\bibitem{KTU} K. Kasamatsu, M. Tsubota and M. Ueda, {\em Giant hole and circular superflow in a fast rotating Bose-Einstein condensate}, Phys. Rev. B {\bf 66} (2002), 053606.
	
	
	

\bibitem{K} M. K. Kwong, {\em Uniqueness of positive solutions of $\Delta u-u+u^p=0$  in $\R^N$}, Arch. Rational Mech. Anal. {\bf 105} (1989), 243--266.
	
	
	
\bibitem{Lewin} M. Lewin, P. T. Nam and N. Rougerie, {\em Blow-up profile of rotating 2D focusing Bose gases}, Macroscopic Limits of Quantum Systems, Springer Verlag, 2018.
	
	
	
	
%\bibitem{LPW} G. Li, S.  Peng and C. Wang, {\em Infinitely many solutions for nonlinear Schrödinger
%		equations with electromagnetic fields}, J. Differential Equations {\bf 251} (2011), 3500--3521.
	
\bibitem{Lieb} E. H. Lieb and M. Loss, Analysis, Graduate Studies in Mathematics Vol. 14, Amer. Math. Soc., Providence, RI, 2001.
	
	
	
\bibitem{Lieb06} E. H. Lieb and R. Seiringer, {\em Derivation of the Gross-Pitaevskii equation for rotating Bose gases}, Comm. Math. Phys. {\bf 264} (2006), 505--537.
	
	
	
	
	
%\bibitem {LSS} E. H. Lieb, R. Seiringer, J. P. Solovej and  J. Yngvason, {\em The mathematics of the Bose gas and its condensation}, Oberwolfach Seminars, {\bf 34} Birkh$\ddot{a}$user Verlag, Basel, 2005.
	
	
	
\bibitem {LSY} E. H. Lieb, R. Seiringer and J. Yngvason, {\em Bosons in a trap: A rigorous derivation of the Gross-Pitaevskii energy functional}, Phys. Rev. A {\bf 61} (2000), 043602.
	
\bibitem {Lions1} P. L. Lions, {\em The concentration-compactness principle in the calculus of variations. The locally compact case}, Part I: Ann. Inst. H. Poincar\'{e} Anal. Non Lin\'{e}aire  {\bf 1} (1984), 109--145. Part II: Ann. Inst. H. Poincar\'{e} Anal. Non Lin\'{e}aire  {\bf 1}  (1984), 223--283.	
	
	%\bibitem {LT} M. Loss and B. Thaller, {\em Optimal heat kernel estimates for Schr\"{o}dinger operators with magnetic field in two dimensions}, Comm. Math. Phys. {\bf 186} (1997), 95--107.
	
\bibitem {LP} K. Lu, X. Pan, {\em Gauge invariant eigenvalue problems in $\R^2$ and in $\R^2_+$}, Trans. Amer. Math. Soc. {\bf 3} (2000), 1247--1276.
	
	
	
	
\bibitem {LC} E. Lundh, A. Collin and K.-A. Suominen, {\em Rotational states of Bose gases with attractive interactions
		in anharmonic traps}, Phys. Rev. Lett. {\bf 92} (2004), 070401.
	
	
%\bibitem {LPY} P. Luo, S. Peng and S. Yan, {\em  Excited states on Bose-Einstein condensates with attractive interactions}, submitted, (2019), arXiv:1909.08828.
%	
%\bibitem{MC1} K. Madison, F. Chevy, J. Dalibard and W. Wohlleben, {\em Vortex formation in a stirred
%		Bose-Einstein condensate}, Phys. Rev. Lett. {\bf 84} (2000), 806.
%	
%\bibitem{MC2} K. Madison, F. Chevy, J. Dalibard and W. Wohlleben, {\em Vortices in a stirred Bose-Einstein condensate}, J. Mod. Opt. {\bf 47} (2000), 2715--2723.
	
	%\bibitem {M} M. Maeda, {\em On the symmetry of the ground states of nonlinear
	
	%Schr$\ddot{o}$dinger equation with potential}, Adv. Nonlinear Stud. {\bf 10}, 895--925 (2010).
	
%\bibitem{NT} W. Ni and I. Takagi, {\em On the shape of least-energy solutions to a semilinear Neumann problem}, Comm. Pure Appl. Math. {\bf 44} (1991), 819--851.
	
	
	
	%\bibitem{NT2} W.-M. Ni and I. Takagi, {\em Locating the peaks of least-energy solutions to a semilinear Neumann problem},  Duke Math. J. {\bf 70} (1993), 247--281.
	%
	%
	%
	%
	%
	%\bibitem{OK} N. Okazawa, {\em An $L^p$ theory for Schr\"{o}dinger operator with nonnegative potentials}, J. Math. Soc. Japan. {\bf 36} (1984), 675--688.
	
	
	
	
	
\bibitem {P} L. P. Pitaevskii, {\em Vortex lines in an imperfect Bose gas}, Sov. Phys. JETP. {\bf 13} (1961), 451--454.
	
	
\bibitem {RS} M. Reed and B. Simon, Methods of modern mathematical physics IV. Analysis of operators, Academic Press, New York-London, 1978.
	%
	%
	%
	%\bibitem{Ro} N. Rougerie, {\em The giant vortex state for a Bose-Einstein condensate in a rotating anharmonic trap: extreme rotation regimes}, J. Math. Pures Appl. {\bf 9} (2011), 296--347.
	
%	
%\bibitem{Hulet3} C. A. Sackett, H. T. C. Stoof and R. G. Hulet, {\em Growth and collapse of a Bose-Einstein condensate with attractive interactions}, Phys. Rev. Lett. {\bf 80} (1998), 2031.
	
	
	
	%\bibitem{SS} E. Sandier and S. Serfaty, {\em Global minimizers for the Ginzburg-Landau functional below the first critical magnetic field}, Ann. Inst. H. Poincar\' e  Anal. Non Lin\' eaire {\bf 17} (2000), 119--145.
	%
	%
\bibitem{SSbook}E. Sandier and S. Serfaty,  Vortices in the Magnetic Ginzburg-Landau Model, Progress in Nonlinear Differential Equations and their Applications {\bf 70}, Basel: Birkh\'auser, 2007.

\bibitem {S02} R. Seiringer, {\em Gross-Pitaevskii theory of the rotating Bose gas}, Comm. Math. Phys. {\bf 229} (2002), 491--509.
	%
	%\bibitem {SS} S. Secchi and M. Squassina, {\em On the location of spikes for the Schr\"{o}dinger equation with electromagnetic field}, Comm. Contemp. Math. {\bf 7} (2005), 251--268.
	
	%
	%%\bibitem {S} R. Seiringer, {\em Hot topics in cold gases}, XVIth International Congress
	%
	%%on Mathematical Physics, World Sci. Publ., Hackensack, NJ, 231--245 (2010).
	%
	%
	%\bibitem {WT} M. Wadati and T. Tsurumi, {\em Critical number of atoms for the magnetically trapped Bose-Einstein condensate with negative s-wave scattering length}, Phys. Lett. A {\bf 247} (1998), 287--293.
	%
	%
	%
	%
	%\bibitem {Wang} X. F. Wang, {\em On concentration of positive bound states of nonlinear Schr$\ddot{o}$dinger equations}, Comm. Math. Phys. {\bf 153} (1993),  229--244.
	
	
	
\bibitem {W} M. I. Weinstein, {\em Nonlinear Schr$\ddot{o}$dinger equations and sharp interpolations estimates}, Comm. Math. Phys. {\bf 87} (1983), 567--576.
	
%\bibitem {WG} N. K. Wilkin, J. M. F. Gunn, and R. A. Smith, {\em Do attractive Bosons condense?}, Phys. Rev. Lett.
%	{\bf 80} (1998), 2265.


\bibitem{WM} M. Willem, Minimax Theorems, Progress in Nonlinear Differential Equations and their Applications 24, Birkh\"{a}user Boston Inc, Boston, (1996).

	
%\bibitem {Zhang} J. Zhang, {\em Stability of standing waves for nonlinear Schr\"{o}dinger equations with unbounded potentials}, Z. Angew. Math. Phys. {\bf 51} (2000), 498--503.
%	
%	
%	
%\bibitem {Z} J. Zhang, {\em Stability of attractive Bose-Einstein condensates}, J. Stat. Phys. {\bf 101} (2000), 731--746.
	
	
	
%\bibitem {ZW} M. W. Zwierlein, J. R. Abo-Shaeer, A. Schirotzek, C. H. Schunck and W. Ketterle, {\em Vortices and superfluidity in a strongly interacting fermi gas}, Nature {\bf 435} (2005), 1047--1051.
	
	
	
\end{thebibliography}
\end{document}